\documentclass[11pt, a4paper]{article}
\usepackage[margin=1.1in]{geometry}
\usepackage[english]{babel}
\usepackage[utf8]{inputenc}
\usepackage[overlay,absolute]{textpos}
\usepackage[T1]{fontenc}
\usepackage{lmodern}
\usepackage{authblk}
\usepackage{url}
\usepackage{lscape}
\usepackage{setspace}
\usepackage{caption}
\usepackage[pdftex]{graphicx}
\usepackage{amssymb,amsmath}
\usepackage{amsthm}
\usepackage[dvipsnames]{xcolor}
\usepackage{hyperref}
\usepackage{mathcomp}
\usepackage{multicol}
\usepackage{ytableau}
\usepackage[nottoc]{tocbibind}
\usepackage{enumitem}
\usepackage{comment}
\usepackage{stmaryrd}
\usepackage{array}
\usepackage{tabularx}
\usepackage{colortbl}
\usepackage{afterpage}
\usepackage{float}
\usepackage{tikz}
\usepackage{tikz-cd}
\usepackage{bm}
\usetikzlibrary{matrix}
\usepackage[all]{xy}
\usepackage{mdframed}
\usepackage{framed}
\usepackage{etoolbox}

\newtheorem{thm}{Theorem}[section]
\newtheorem{prop}[thm]{Proposition}
\newtheorem{cor}[thm]{Corollary}
\newtheorem{lem}[thm]{Lemma}

\newcommand{\A}{\mathcal{A}}
\newcommand{\St}{\mathcal{S}}
\newcommand{\B}{\mathcal{B}}
\newcommand{\dom}{\mathrm{dom}\,}
\newcommand{\im}{\mathrm{im}\,}

\newcommand{\rr}{\mathrm}
\newcommand{\cl}{\mathcal}
\newcommand{\br}[1]{\boldsymbol{\mathbf{#1}}}

\mdfsetup{skipabove=2pt,skipbelow=8pt}
\newmdtheoremenv[linewidth=0.9pt, topline=false, bottomline=false, rightline=false,%
leftmargin=0pt, innerleftmargin=0.4em, rightmargin=0pt, innerrightmargin=0pt, innertopmargin=-5pt ,%
innerbottommargin=3pt, splittopskip=\topskip, splitbottomskip=0.3\topskip, %
skipabove=0.6\topsep]%
{exmp}[thm]{Example}%

\newcommand{\red}{\textcolor{Red}}
\newcommand{\blue}{\textcolor{NavyBlue}}
\definecolor{myblue}{RGB}{24,65,223}
\definecolor{myred}{RGB}{229,21,1}
\definecolor{mypurple}{RGB}{126,43,112}

\textblockorigin{-5mm}{0mm} 

\makeatletter
\newcommand*{\textoverline}[1]{$\overline{\hbox{#1}}\m@th$}
\makeatother

\makeatletter
\renewcommand\subitem{\@idxitem\nobreak\hspace*{20\p@}}
\makeatother

\makeatletter \g@addto@macro\@floatboxreset\centering \makeatother

\begin{document}

\title{On the transition monoid of the Stallings automaton of a subgroup of a free group}
\author{\textbf{Inês F. Guimarães}}
\affil{\textit{Centro de Matem\'{a}tica, Faculdade de Ci\^{e}ncias, Universidade do
Porto,} \\ \textit{R. Campo Alegre 687, 4169-007 Porto, Portugal} \\  inesguimaraes42@gmail.com}
\maketitle


\section*{Abstract}

Birget, Margolis, Meakin and Weil proved that a finitely generated subgroup $K$ of a free group is pure if and only if the transition monoid $M(K)$ of its Stallings automaton is aperiodic. In this paper, we establish further connections between algebraic properties of $K$ and algebraic properties of $M(K)$. We mainly focus on the cases where $M(K)$ belongs to the pseudovariety of finite monoids all of whose subgroups lie in a given pseudovariety of finite groups. We also discuss normal, malnormal and cyclonormal subgroups of $F_A$ using the transition monoid of the corresponding Stallings automaton.

\vspace{2mm}

\noindent \textbf{Keywords:} free group, Stallings automaton, transition monoid, pseudovariety, normal subgroup, malnormal subgroup

\vspace{2mm}

\noindent 2010 Mathematics Subject Classification: 20E07, 20M07, 20M35, 20E36, 68Q45

\section{Introduction}

With the purpose of finding efficient methods to tackle problems involving subgroups of free groups, John Stallings \cite{stallings} presented in a paper of 1983 a revolutionary approach. He developed a way of associating with each finitely generated subgroup of a free group a finite labeled graph, under the formalism of graph immersions. These graphs became known as \textbf{Stallings automata} and constitute a most powerful tool for studying finitely generated subgroups (``f.g. subgroups'' for short) of a free group. 

Some years later, Alexei Miasnikov and Ilya Kapovich gave Stallings' construction a more combinato\-rial flavor and collected numerous group-theoretic properties of f.g. subgroups of free groups based on combinatorial properties of their Stallings automata \cite{kapovich}. Besides being a very elegant theory, this approach displayed great benefits from an algorithmic viewpoint. For a list of applications of Stallings automata, see \cite{list}. 

Furthermore, the paper \cite{inspo} written by Birget, Margolis, Meakin and Weil unveiled a new way of characterizing properties of a finitely generated subgroup $K$ of a free group, this time by looking at algebraic properties of its Stallings automaton $\St(K)$, namely its transition monoid $M(K)$. In this article, the authors prove that the condition of a f.g. subgroup $K$ of a free group being pure (respectively, $p$-pure) is equivalent to the condition of $M(K)$ being aperiodic (respectively, $p$-periodic). Clearly, there are important instances of pseudovarieties of monoids involved in this result. 

Indeed, finite automata and finite monoids are deeply related and constitute important objects not only in mathematics but also in computer science. Moreover, pseudovarieties and varieties contribute greatly to the classification of finite monoids and rational lan\-guages. Hence, given a f.g. subgroup $K$ of a free group, there is a strong motivation to find further pseudovarietal properties of $M(K)$ (other than being aperiodic or $p$-periodic) that correspond to algebraic properties of $K$ (as a subgroup of the ambient free group). More concretely, let $A$ be a finite alphabet and let $K \leq_{f.g.} F_A$ be a finitely generated subgroup of the free group $F_A$ over $A$. For what pseudovarieties of monoids $\br{V}$ can we guarantee that $M(K) \in \br{V}$ implies $M(K\varphi) \in \br{V}$ for every automorphism $\varphi$ of $F_A$? We say that such pseudovarieties are \textbf{stable under every automorphism of $F_A$}. And what do those ``well-behaved'' pseudovarieties tell us about $K$ itself? 

Our main results concern pseudovarieties of a certain type that satisfy the previous property and some characterizations of the subgroups involved. Given a pseudovariety of groups $\boldsymbol{\mathrm{H}}$, we denote by $\overline{\boldsymbol{\mathrm{H}}}$ the pseudovariety of monoids all of whose subgroups lie in $\boldsymbol{\mathrm{H}}$. In Theorem~\ref{thm:thebigone}, we prove that if $\boldsymbol{\mathrm{H}}$ is a pseudovariety of groups, then the pseudovariety of monoids $\overline{\boldsymbol{\mathrm{H}}}$ is stable under all automorphisms of $F_A$. Observe that the pseudovariety of aperiodic monoids is precisely the pseudovariety of monoids all of whose subgroups are trivial, that is, $\boldsymbol{\mathrm{A}} = \overline{\boldsymbol{\mathrm{{I}}}}$.
Let $S = (k_n)_{n \geq 1}$ be a sequence of positive integers and denote by $\boldsymbol{\mathrm{V}}_S$ the pseudovariety of finite groups ultimately defined by the sequence $(x^{k_n}=1)_n$. In Theorem~\ref{thm:psv:xkn}, it is shown that $M(K) \in \overline{\boldsymbol{\mathrm{V}}}_S$ if and only if there exists  $p\geq 1$ such that for all $x \in F_A$, $n \geq 1$ and $m \geq p$ we have $x^n \in K \Rightarrow x^{(n,k_m)} \in K$.

Furthermore, inspired by the characterizations of normal, malnormal and cyclonormal subgroups in terms of their Stallings automata presented in \cite{kapovich}, we were also motivated to discuss these conjugacy conditions by inspecting the structure of the transition monoids of Stallings automata. More precisely, let $Q$ be the vertex set of the automaton $\mathcal{S}(K)$ and, given $u \in (A \cup A^{-1})^*$, denote by $\delta_u$ the partial transformation over $Q$ mapping a vertex $q \in Q$ to the vertex reached after reading $u$ from $q$ in $\mathcal{S}(K)$, whenever that is possible. In Theorem~\ref{thm:size:group}, we prove that a nontrivial subgroup $K \leq_{f.g.} F_A$ is normal in $F_A$ if and only if $M(K)$ is a group of size $|Q|$. Moreover, let $R_A$ be the set of all reduced words over the alphabet $A \cup A^{-1}$, and let $E = \{\delta_u \in M(K) \mid \delta_u = \delta_u^2, \, u \in R_A \setminus \{1\} \}$. Consider the restriction of the natural partial order on $M(K)$ to $E$, and let $k$ be the size of a maximal chain on $E$. In Theorem~\ref{thm:malnormal:me}, we prove that a nontrivial subgroup $K <_{f.g.} F_A$ is malnormal if and only if $k=2$ and $|E| = |Q| + 1$.

This paper is organized as follows. In Section~\ref{sec:preliminaries}, we present some background concepts needed for the upcoming sections. In Section~\ref{sec:auto}, we include some results regarding the structure of the automorphism group of a free group, as well as the effect of automorphisms at the level of Stallings automata. We briefly explore what happens when the transition monoid of a Stallings automaton is a group in Section~\ref{sec:transgroup}. In Section~\ref{sec:main}, we state our most important results regarding pseudovarieties that are stable under every automorphism of a free group. Finally, in Section~\ref{sec:conjug}, we discuss normal, malnormal and cyclonormal subgroups of $F_A$ using the transition monoid of the Stallings automaton. We end with some comments regarding future directions of work.

\section{Preliminaries}\label{sec:preliminaries}

\subsection{Free groups}

Let $A$ be a finite alphabet. We denote by $A^*$ the \textbf{free monoid} over $A$, by $A^{-1}$ the set of formal inverses of the letters in $A$, and we write $\widetilde{A} = A \cup A^{-1}$.

By successively erasing factors of the type $aa^{-1} (a \in \widetilde{A})$ from a word $w \in \widetilde{A}^*$, we arrive at the unique \textbf{reduced} word $\overline{w}$ without factors of this kind. We can then consider the congruence $\tau_A \subseteq \widetilde{A}^* \times \widetilde{A}^*$ given by \[(u,v) \in \tau_A \Leftrightarrow \overline{u} = \overline{v}.\] Finally, we define $F_A = \{u \tau_A \mid u \in \widetilde{A}^* \}$ which, when endowed with the binary operation $(u \tau_A)(v \tau_A) = (uv) \tau_A$, becomes the \textbf{free group} over $A$.

When no confusion arises, we write $\overline{u}$ or even $u$ instead of $u \tau_A$; in particular, we consider $\widetilde{A}$ as a subset of $F_A$. Moreover, given a subgroup $K \leq F_A$, we write $\overline{K}$ to designate the set of all reduced words in $\widetilde{A}^*$ representing the elements of $K$. 

\subsection{Finite automata and rational languages}

Let $A$ be a finite alphabet. A \textbf{language} over $A$ (or $A$-language) is a subset of $A^*$. We can combine $A$-languages using the so-called rational operators: union, product and star. The star operator applied to a language $L$ yields the language \[ L^* = \{u_1\dots u_n \mid n \geq 0; u_1, \dots, u_n \in L\} = \bigcup_{n \geq 0} L^n,\] under the convention $L^0 = \{1\}$, where $1$ denotes the empty word.

Consider the congruence $\sim_L$ on $A^*$ given by 
\[u \sim_L v \; \text{ if } \; \forall x,y \in A^* \; (xuy \in L \Leftrightarrow xvy \in L).\]
We define the \textbf{syntactic monoid} of $L$ to be $\mathrm{Synt}(L) = A^*/{\sim_L}$.

An $A$-language is called \textbf{rational} if it can be obtained from finite $A$-languages using the rational operators finitely many times, which is the same as saying that it admits a rational expression.

A finite \textbf{automaton} over $A$, also called an $A$-automaton, is a structure of the type $\A = (Q,A,E,I,T)$, where:
\begin{itemize}\setlength\itemsep{-0.1em}
\item $Q$ is a finite set, called the set of \textbf{vertices} or \textbf{states};
\item $E \subseteq Q \times A \times Q$ is the (finite) set of \textbf{edges} or \textbf{transitions};
\item $I, T \subseteq Q$ are the sets of \textbf{initial} and \textbf{terminal} states, respectively.
\end{itemize}

When it is possible to read every letter of $A$ from any vertex in $Q$, we say that $\A$ is a \textbf{complete} automaton. The \textbf{underlying graph} of $\mathcal{A}$ is the directed labeled graph obtained from $\A$ by ignoring the designation of vertices as initial or terminal; we denote it by $\Gamma_{\A}$.

A \textbf{path} in an automaton $\A = (Q,A,E,I,T)$ is a sequence of the type
\[(p_0,a_1,p_1)(p_1,a_2,p_2)\cdots(p_{n-1},a_n,p_n)\]
where $n\geq 0$, $p_i \in Q$ ($0 \leq i \leq n$) and $(p_{i-1}, a_i, p_i) \in E$ ($1 \leq i \leq n$). We represent it by
\[\xygraph{  
!{<0cm,0cm>;<1cm,0cm>:<0cm,1cm>::}  
!{(0,0) }*+{p_0}="p0"  
!{(1.5,0) }*+{p_1}="p1"   
!{(3,0) }*+{\cdots}="..."
!{(4.5,0) }*+{p_n}="pn"
"p0":"p1"^{a_1}
"p1":"..."^{a_2} 
"...":"pn"^{a_n}
}.\] 
We call $p_1, \dots, p_{n-1}$ the intermediate vertices of the path. If $n>0$, the \textbf{label} of such a path is the word $u = a_1a_2\dots a_n \in A^*$; if $n=0$, we get the trivial path at $p_0$ and its label is the empty word $1 \in A^*$. If there exists a path between vertices $p,q \in Q$ labeled by $u \in A^*$, we represent it by $\xymatrix{p \ar[r]^u & q}$; and when $p=q$, we say that $u$ labels a \textbf{loop} at $p$.

A path $p_0 \xrightarrow{\; a_1 \;} p_1 \xrightarrow{\; a_2 \;} \cdots \xrightarrow{\; a_n \;} p_n$ in $\A$ is called \textbf{successful} if $p_0 \in I$ and $p_n \in T$. We define the \textbf{language recognized by $\A$} as
\[L(\A) = \{u \in A^* \mid u \text{ is the label of a successful path in } \A\}. \] 
The following result is well known (for a proof, see \cite[Chapter 5, Theorem 5.2.1]{lawson}) and establishes a connection between finite automata and rational languages.

\begin{thm}[Kleene's Theorem]\label{thm:kleene}
Let $A$ be a finite alphabet and $L \subseteq A^*$. Then there exists a finite automaton $\A$ satisfying $L(\A) = L$ if and only if $L$ is a  rational language.
\end{thm}

We call an automaton $\A = (Q,A,E,I,T)$ \textbf{deterministic} if it has a unique initial state and
\[ (p,a,q), (p,a,r) \in E \Rightarrow q = r\]
for all $p,q,r \in Q$ and $a \in A$. When $\A$ is deterministic, we define a partial function $\delta$ on $Q \times A$ by
\[ (p,a)\delta = q \Leftrightarrow (p,a,q) \in E \]
for all $p,q \in Q$ and $a \in A$. We call it the \textbf{transition function} of $\A$ and we often write $\A = (Q,A,\delta,q_0,T)$ instead of $\A = (Q,A,E,q_0,T)$ when dealing with deterministic automata. We denote by $0$ the empty transformation, that is, the partial transformation whose domain is the empty set $\emptyset$.

We can extend $\delta$ to $Q \times A^*$ by letting $(q,u)\delta$ be the state reached after reading the word $u \in A^*$ from $q \in Q$ by following the labels on the edges (whenever that is possible). When no confusion arises, we write $q \cdot u$ instead of $(q,u)\delta$.

Given $u \in A^*$, we further define $\delta_u \colon Q \rightarrow Q$ by $q \delta_u = q \cdot u$. Denoting by $\mathcal{PT}_Q$ the monoid of all partial transformations on $Q$, it is easy to check that the map
\begin{align*}
\Delta \colon A^* &\longrightarrow \mathcal{PT}_Q \\
u &\longmapsto \delta_u
\end{align*}
is a monoid homomorphism; hence, the image of $\Delta$ is a submonoid of $\mathcal{PT}_Q$. We call it the \textbf{transition monoid} of $\A$ and denote it by $M(\A)$.

Given an automaton $\A = (Q,A,E,I,T)$, a subset $P \subseteq Q$ and a word $u \in A^*$, we define 
\[Pu = \{q \in Q \mid \text{there exists a path} \xymatrix{p \ar[r]^u & q}\, \text{in } \A \text{ for some } p \in P\}. \]
A vertex $q \in Q$ is called \textbf{accessible} if there exists $u \in A^*$ such that $q \in Iu$ and \textbf{co-accessible} if $qu \in T$ for some $u \in A^*$. Clearly, eliminating vertices which are not acessible and vertices which are not co-accessible does not change the language recognized by the automaton. When all vertices are both accessible and co-accessible, the automaton is said to be \textbf{trim}. 

An automaton $\A = (Q,\widetilde{A},E,I,T)$ is called \textbf{involutive} if, for all $p,q \in Q$ and $a \in A$, we have
\[(p,a,q) \in E \Leftrightarrow (q,a^{-1},p) \in E.\]
These pairs of edges are considered inverses of each other, and an edge labeled by a letter in $A$ (respectively, in $A^{-1}$) is said to be positive (respectively, negative). When we know beforehand that a certain automaton is involutive, we only draw the positive edges; the negative ones are like ``ghost'' edges that we visualize only in our heads. Moreover, we denote by $E^+ \subseteq E$ and $E^- \subseteq E$ the subsets of positive and negative edges, in that order. 

Let $\A$ be an involutive automaton over $A$ and let $w = a_1 a_2 \dots a_k \in \widetilde{A}^*$, with $a_i \in \widetilde{A}$ for $1 \leq i \leq k$. Suppose that $w$ labels a path 
\[\xymatrix{p_0 \ar[r]^{a_1} & p_1 \ar[r]^{a_2} & \dots \ar[r]^{a_k} & p_k}\]
in $\A$. We say that the path labeled by $w$ is \textbf{reduced} if it does not contain any consecutive edges of the form
\[\xymatrix{p_{i-1} \ar[r]^{a_i} & p_i \ar[r]^{a_i^{-1}} & p_{i-1}}\]
for $1 \leq i \leq k$.

An \textbf{inverse} automaton is an involutive, deterministic and trim automaton with a unique terminal state. It is folklore that the transition monoid $M(\A)$ of an inverse automaton $\A$ is an \textbf{inverse monoid}, i.e., for every $x \in M(\A)$, there exists a unique $x^{-1} \in M(\A)$ satisfying $xx^{-1}x = x$ and $x^{-1}x x^{-1} = x^{-1}$. Moreover, since an inverse automaton is a minimal automaton \cite{bartholdi}, it follows that the transition monoid of an inverse automaton $\A$ is isomorphic to the syntactic monoid of the language recognized by $\A$, that is, $M(\A) \cong \mathrm{Synt}(L(\A))$ \cite{pin}. For more details regarding rational languages and automata, the reader is referred to \cite{eilenbergA, pin}.

Inverse monoids, being the transition monoids of inverse automata, will play an important role in what follows, so we end this section by presenting some equivalence relations $\cl{R}$, $\cl{L}$, $\cl{H}$ and $\cl{D}$, known as \textbf{Green's relations} \cite{green}, which are useful to analyze the structure of an inverse monoid.

Let $M$ be an inverse monoid and $x,y \in M$. We say that $x$ and $y$ are: 
\begin{itemize}\setlength\itemsep{-0.1em}
\item $\mathcal{R}$-related if $xx^{-1} = yy^{-1}$;
\item $\mathcal{L}$-related if $x^{-1}x = y^{-1}y$;
\item $\mathcal{H}$-related if $xx^{-1} = yy^{-1}$ and $x^{-1}x = y^{-1}y$;
\item $\mathcal{D}$-related if there exists $z \in M$ such that $x\mathrel{\mathcal R}z$ and $z\mathrel{\mathcal L}y$.
\end{itemize}
If $\cl{K}$ is one of Green's relations, we write $x \, \cl{K} \, y$ to indicate that $x$ and $y$ are $\cl{K}$-related, i.e., belong to the same $\cl{K}$-class. A monoid $M$ is called $\cl{K}$-trivial if all of its $\cl{K}$-classes are singletons. For more details regarding Green's relations and the structure of inverse monoids, the reader is referred to \cite{howie, petrich}.

\subsection{Stallings' construction}

Let $F_A$ be a free group over a finite alphabet $A$. Let $K = \langle u_1, \dots, u_k \rangle \leq F_A$ be a finitely generated subgroup of $F_A$, where each generator $u_i$ is seen as a (nonempty) reduced word in $\widetilde{A}^*$. The notation $\leq_{f.g.}$ will be often used to indicate that a subgroup is finitely generated. We begin by constructing the so-called \textbf{flower automaton} $\mathcal{F}(u_1,\dots,u_k)$ of $K$ by fixing a point $q_0$, called the \textbf{basepoint} of the automaton, and gluing to it $k$ ``petals'' labeled by each of the $u_i$, as well as the corresponding inverse edges, in order to obtain an involutive automaton over $\widetilde{A}$:
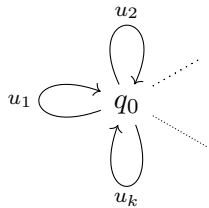
\vspace{2mm}
\[
\begin{tikzcd}
    q_0\ar[
    ,loop 
    ,out=115 
    ,in=65 
    ,distance=3em 
    ]{}{u_2}\ar[
    ,loop 
    ,out=-70 
    ,in=-110 
    ,distance=3em 
    ]{}{u_k}\ar[
    ,loop 
    ,out=-160 
    ,in=-200 
    ,distance=3em 
    ]{}{u_1}\ar[out=-30, in=-30, distance = 3em, no head, dotted]\ar[out=30, in=30, distance = 3em, no head, dotted] & 
\end{tikzcd}
\]
We declare $q_0$ to be the unique initial state and the unique terminal state of the flower automaton. To turn this into an inverse automaton, whenever we encounter a pair of distinct edges $\xymatrix{p \ar[r]^-{a} & q}$ and $\xymatrix{p \ar[r]^-{a} & r}$ for some $a \in \widetilde{A}$, we identify them, and we also identify the corresponding inverse edges (so $q$ and $r$ collapse into a single vertex if they are distinct). These identifications are known as \textbf{Stallings foldings}, and they are successively applied until we reach a deterministic automaton. The inverse automaton thus obtained is called a \textbf{Stallings automaton} of $K$.

\begin{prop}\label{prop:stallingslanguage}
Let $F_A$ be a free group of finite rank and let $K \leq_{f.g.} F_A$. Then the language recognized by any Stallings automaton of $K$ is the intersection of all languages $L \subseteq \widetilde{A}^*$ containing $\overline{K}$ which are recognized by a finite inverse automaton with a basepoint.
\end{prop}

This result is proven in \cite{bartholdi}. It allows us to conclude that any two Stallings automata of $K$ are isomorphic, so we can speak of \textit{the} Stallings automaton of $K$, denoting it by $\St(K)$. In other words, $\St(K)$ does not depend on the generating set of $K$ nor on the order in which the foldings are made. However, $\St(K)$ depends on the basis $A$ of the free group we are considering.

We illustrate Stallings' construction with an example.

\begin{exmp}\label{ex:stallings1}
\normalfont Let $A = \{a,b,c\}$ and $K = \langle \textcolor{Red}{c}, \textcolor{Green}{ba^{-1}c^{-1}}, \textcolor{NavyBlue}{aca^{-1}} \rangle \leq F_A$. Then the flower automaton $\cl{F}(c,ba^{-1}c^{-1},aca^{-1})$ is depicted by 
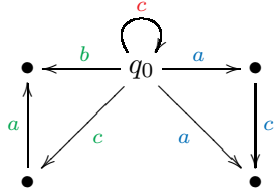
\[\xygraph{  
!{<0cm,0cm>;<1cm,0cm>:<0cm,1cm>::}  
!{(0,0)}*+{q_0}="0"  
!{(-1.5,0)}*+{\bullet}="-1"   
!{(1.5,0)}*+{\bullet}="1" 
!{(-1.5,-1.5)}*+{\bullet}="-2"
!{(1.5,-1.5)}*+{\bullet}="2"
"0":@[NavyBlue]"1"^{\textcolor{NavyBlue}{a}} "1":@[NavyBlue]"2"^{\textcolor{NavyBlue}{c}} "0":@[NavyBlue]"2"_{\textcolor{NavyBlue}{a}}
"0":@[Green]"-1"_{\textcolor{Green}{b}} "0":@[Green]"-2"^{\textcolor{Green}{c}} "-2":@[Green]"-1"^{\textcolor{Green}{a}}
"0" :@(lu,ru)@[Red]"0"^{\textcolor{Red}{c}}
} \]
After folding the two blue edges labeled by $a$ and the green and red edges labeled by $c$, we obtain the automaton
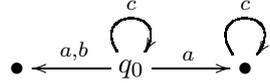
\[\xygraph{  
!{<0cm,0cm>;<1cm,0cm>:<0cm,1cm>::}  
!{(0,0)}*+{q_0}="0"  
!{(-1.5,0)}*+{\bullet}="-1"   
!{(1.5,0)}*+{\bullet}="1" 
"0":"1"^{a}
"0":"-1"_{a,b}
"0":@(lu,ru)"0"^{c}
"1":@(lu,ru)"1"^{c}
} \]
It remains to fold the edges labeled by $a$, so $\St(K)$ is given by
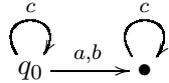
\[\xygraph{  
!{<0cm,0cm>;<1cm,0cm>:<0cm,1cm>::}  
!{(0,0)}*+{q_0}="0"    
!{(1.5,0)}*+{\bullet}="1" 
"0":"1"^{a,b}
"0":@(lu,ru)"0"^{c}
"1":@(lu,ru)"1"^{c}
} \]
\end{exmp}

The next result follows from the proof of Proposition~\ref{prop:stallingslanguage}, and it allows us to conclude that the \textbf{generalized word problem} for finitely generated free groups is decidable (for details, see \cite{bartholdi}).

\begin{prop}\label{prop:stallingslanguage2}
Let $K \leq_{f.g.} F_A$. Then $\overline{L(\St(K))} = \overline{K}$ and $u \in F_A$ belongs to $K$ if and only if $\overline{u} \in L(\St(K))$.
\end{prop}

\subsection{Pseudovarieties}\label{subsec:psv}

A \textbf{pseudovariety of (finite) monoids} is a class of (finite) monoids $\br{V}$ closed under taking submonoids, homomor\-phic images and (finitary) direct products. This means that:
\begin{enumerate}[label=(\roman*)]\setlength\itemsep{-0.1em}
\item For all $M \in \br{V}$, if $N \leq M$, then $N \in \br{V}$.
\item For all $M \in \br{V}$, if $\varphi \colon M \rightarrow N$ is an onto monoid homomorphism, then $N \in \br{V}$.
\item For all $M,N \in \br{V}$, $M \times N \in \br{V}$.
\end{enumerate}

There is yet another way to describe pseudovarieties of monoids, according to the ``equa\-tions'' they satisfy. Given a set of variables $A$, a \textbf{monoid identity} on $A$ is an element $(u,v) \in A^* \times A^*$, which we usually indicate by a formal equality $u = v$. We say that a monoid $M$ satisfies the identity $u = v$ if $u\varphi = v\varphi$ for every monoid homomorphism $\varphi \colon A^* \rightarrow M$. In that case, we write $M \models u=v$. Informally, this means that we obtain a true equality when we replace the variables of $A$ by arbitrary elements of $M$ in the identity $u=v$. If $(u_n = v_n)_n$ is a sequence of monoid identities, we say that $M$ \textbf{ultimately satisfies} that sequence of identities if there exists some $p\geq 1$ such that $M \models u_n = v_n$ for every $n \geq p$. If $\cal{C}$ is the class of monoids ultimately satisfying the sequence of identities $(u_n = v_n)_n$, we say that $\cal{C}$ is \textbf{ultimately defined} by $(u_n = v_n)_n$. The next result provides a characterization of pseudovarieties of monoids in terms of monoid identities, and a proof can be found in \cite{schutzenberger}.

\begin{thm}[Eilenberg and Schützenberger]\label{thm:psvmidentities}
A class of finite monoids $\boldsymbol{\mathrm{V}}$ is a pseudovari\-ety of monoids if and only if $\boldsymbol{\mathrm{V}}$ is ultimately defined by a sequence of identities $(u_n = v_n)_n$. 
\end{thm}

We now introduce a notation that will be useful to simplify the writing of monoid identities. Given a monoid $M$, we denote by $E(M)$ the set of \textbf{idempotents} of $M$, i.e. $E(M) = \{ x \in M \mid x^2 = x\}$. Observe that in a finite monoid there exists $k \geq 1$ such that $x^k \in E(M)$ for all $x \in M$; we call $k$ an \textbf{exponent} of $M$. If $\overline{n}$ denotes the least common multiple of the numbers $1,2,\dots,n$, then $x^{\overline{n}}$ is the unique idempotent power of $x$ whenever $n \geq k$, for some exponent $k$ of $M$. Following a convention of Sch\"utzenberger, in a pseudovariety defined by the identities $(u_n = v_n)_n$, we agree to replace every occurrence of $\overline{n}$ by the symbol $\omega$. If $\Sigma$ is a set of monoid identities (which may feature the $\omega$ symbol), we denote by $\llbracket \Sigma \rrbracket$ the pseudovariety of monoids satisfying all the identities in $\Sigma$.

\begin{exmp}\label{ex:stallings2}
\normalfont Let $\br{G}$, $\br{Com}$, $\br{Sl}$ and $\br{A}$ denote the pseudovarieties of all finite groups, commutative monoids, semilattices and aperiodic monoids, respectively. Then:

\begin{enumerate}[label=(\roman*)]\setlength\itemsep{-0.1em}

\item $\br{G} = \llbracket x^\omega = 1 \rrbracket$

\item $\br{Com} = \llbracket xy = yx \rrbracket$

\item $\br{Sl} = \llbracket x^2 = x, xy = yx \rrbracket $

\item $\br{A} = \llbracket x^{\omega + 1} = x^\omega \rrbracket$

\end{enumerate}

\end{exmp}

As in the case of monoids, a \textbf{pseudovariety of finite groups} is a class of (finite) groups closed under taking subgroups, homomor\-phic images and (finitary) direct products. Equivalently, it is a pseudovariety of monoids whose elements are groups.

Given a pseudovariety of groups $\br{H}$, we denote by $\overline{\br{H}}$ the pseudovariety of finite monoids all of whose subgroups lie in the pseudovari\-ety of groups $\boldsymbol{\mathrm{H}}$. As the group $\mathcal{H}$-classes of a monoid are precisely its maximal subgroups and pseudovarieties of groups are closed under taking subgroups, we can also say that $\overline{\boldsymbol{\mathrm{H}}}$ is the pseudovariety of monoids all of whose group $\mathcal{H}$-classes belong to $\boldsymbol{\mathrm{H}}$. Pseudovarieties of this kind will play a major role in this paper. For more details concerning pseudovarieties of monoids, the reader is refered to \cite{jalmeida, pin}.

\section{Free group automorphisms}\label{sec:auto}

We now discuss the automorphism group of a free group, that is, the group $\rr{Aut}(F_A)$ whose elements are the automorphisms of $F_A$, for some finite alphabet $A$. Observe that an automorphism in $\rr{Aut}(F_A)$ maps any basis of $F_A$ to another basis of $F_A$, and it is completely determined by the images of the elements of a basis. The next result exhibits a finite generating set of $\rr{Aut}(F_A)$.

\begin{thm}\label{thm:nielsengenerators}
Let $A$ be a finite alphabet with at least two elements and $F_A$ the free group over $A$. The automorphisms of the form

\vspace{-6mm}

\hspace{15mm}\begin{minipage}[t]{.4\textwidth}
\begin{align*}
\alpha_{a} \colon F_A &\longrightarrow F_A \\
a &\longmapsto a^{-1} \\
x &\longmapsto x \quad (x \in A \setminus\{a\}) 
\end{align*}
\end{minipage}
\begin{minipage}[t]{.4\textwidth}
\begin{align*}
\beta_{ab} \colon F_A &\longrightarrow F_A \\
a &\longmapsto ab \\
x &\longmapsto x \quad (x \in A \setminus\{a\}),  
\end{align*}
\end{minipage}

\vspace{5mm}

\noindent with $a,b \in A$ distinct letters, generate $\mathrm{Aut}(F_A)$ as a group.
\end{thm}

A proof can be found in \cite[Proposition 4.1]{lyndon}, and it was Nielsen \cite{nielsen} who first presented a set of generators of $\rr{Aut}(F_A)$ very similar to the one above. For that reason, we designate such automorphisms by \textbf{elementary Nielsen automorphisms}. Moreover, we call $\alpha_a$ an automorphism of \textbf{type 1} and $\beta_{ab}$, $\beta_{ab}^{-1}$ automorphisms of \textbf{type 2}. Note that $\beta_{ab}^{-1}$ is the automorphism whose restriction to $A$ consists of replacing $a$ by $ab^{-1}$ and fixing all other letters of $A$. Such a manageable generating set of $\rr{Aut}(F_A)$ will be absolutely crucial to derive our most important results.

We now investigate the effect of applying an automorphism of $F_A$ to a f.g. subgroup $K \leq_{f.g.} F_A$ at the level of its Stallings automaton $\St(K)$. A fact that will sometimes be useful is that any automorphism $\varphi \in \mathrm{Aut}(F_A)$ induces a mapping $A \rightarrow \widetilde{A}^*$, $a \mapsto \overline{a\varphi}$, which can be uniquely extended to a free monoid endomorphism $\phi \colon \widetilde{A}^* \rightarrow \widetilde{A}^*$. 

\begin{prop}\label{prop:edge:replace}
Let $\mathcal{S}(K) = (Q, A, \delta, q_0)$ be the Stallings automaton of $K \leq_{f.g.} F_A$ and let $\varphi \in \mathrm{Aut}(F_A)$. Denote by $\phi \colon \widetilde{A}^* \rightarrow \widetilde{A}^*$ the monoid homomorphism induced by $\varphi$. Then $\mathcal{S}(K\varphi)$ is (isomorphic to) the inverse automaton obtained from $\mathcal{S}(K)$ by the following procedure:
\begin{enumerate}\setlength\itemsep{-0.1em}
\item For every edge $ q \xrightarrow{\; a \;} q' $ in $\mathcal{S}(K)$, let $a\phi = a_1 a_2 \dots a_k$ ($k \geq 1$, $a_i \in \widetilde{A}$) be the factorization of the word $a\phi$ into letters. Then replace that edge by the sequence of edges
\[ q = p_0 \xrightarrow{\,a_1\,} p_1 \xrightarrow{\,a_2\,}  \dots \xrightarrow{\,a_k\,} p_k = q', \]
where $p_1,\dots,p_{k-1}$ are new vertices, and add the corresponding inverse edges.
\item After completing all the edge replacements mentioned above, apply the necessary fol\-dings in order to get an inverse automaton. 
\item Successively eliminate every vertex with outdegree $1$ which is not the basepoint.
\end{enumerate}
\end{prop}

\begin{proof}
Let $\mathcal{A}$ be the automaton we obtain after following step $1$; $\mathcal{A}'$ the automaton we get by the end of step $2$; and $\mathcal{A}''$ the final automaton.

Since $\mathcal{A}''$ is inverse and the only vertex which may have outdegree $1$ is the basepoint, we know that it is the Stallings automaton of some finitely generated subgroup of $F_A$. Thus, we only need to show that $\overline{K\varphi} = \overline{L(\mathcal{A}'')}$, in view of Proposition~\ref{prop:stallingslanguage2}. We start by observing that $\overline{L(\mathcal{A}'')} = \overline{L(\mathcal{A}')} = \overline{L(\mathcal{A})}$. Indeed, regarding the latter equality, given $u \in L(\mathcal{A}')$, we know that $u$ labels a successful path in $\mathcal{A}'$ which can be lifted to a successful path in $\mathcal{A}$ labeled by a word $v$ obtained by inserting factors of the form $aa^{-1}$ ($a \in \widetilde{A}$) into $u$. Hence, $\overline{u} = \overline{v} \in \overline{L(\mathcal{A})}$ and $\overline{L(\mathcal{A}')} \subseteq \overline{L(\mathcal{A})}$. On the other hand, since the inverse automaton $\mathcal{A}'$ is obtained from $\mathcal{A}$ by simply folding edges, it is clear that $L(\mathcal{A}) \subseteq L(\mathcal{A}')$, so $\overline{L(\mathcal{A})} \subseteq \overline{L(\mathcal{A}')}$. As for the first equality, since $\mathcal{A}''$ is a subautomaton of $\mathcal{A}'$, we have $L(\mathcal{A}'') \subseteq L(\mathcal{A}')$, and any word accepted by $\mathcal{A}'$ gives rise to a word accepted by $\mathcal{A}''$ by deleting some factors of the form $aa^{-1}$ ($a \in \widetilde{A}$). Therefore, $\overline{L(\mathcal{A}')} \subseteq \overline{L(\mathcal{A}'')}$ and the chain of equalities follows.

Now, given $u \in \overline{K\varphi}$, there exists $v = v_1 v_2 \dots v_k \in \overline{K}$, with $v_i \in \widetilde{A}$ for $1 \leq i \leq k$, such that $u = \overline{v\phi}$. We know that $v$ labels a path \[ q_0 \xrightarrow{\,v_1\,} q_1 \xrightarrow{\,v_2\,}  \cdots \xrightarrow{\,v_k\,} q_k = q_0 \] in $\mathcal{S}(K)$ so, by step $1$, we get a successful path \[ q_0 \xrightarrow{\,v_1\phi\,} q_1 \xrightarrow{\,v_2\phi\,}  \cdots \xrightarrow{\,v_k\phi\,} q_k = q_0 \] in $\mathcal{A}$. We deduce that $v\phi \in L(\mathcal{A})$, which entails $u = \overline{v\phi} \in \overline{L(\mathcal{A})}$. Hence, $\overline{K\varphi} \subseteq \overline{L(\mathcal{A})}$.

As for the opposite inclusion, given $u = u_1 u_2 \dots u_k \in L(\mathcal{A})$, with $u_i \in \widetilde{A}$ for $1 \leq i \leq k$, we have a path \[ q_0 \xrightarrow{\,u_1\,} q_1 \xrightarrow{\,u_2\,}  \cdots \xrightarrow{\,u_k\,} q_k = q_0 \] in $\mathcal{A}$, which can be modified to yield a reduced path $q_0 \xrightarrow{\;w\;} q_0$ in $\mathcal{A}$. By doing that, we get $\overline{u} = \overline{w}$, and there exists a successful path $q_0 \xrightarrow{\;v\;} q_0$ in $\mathcal{S}(K)$ for some $v \in \widetilde{A}^*$ satisfying $w = v\phi$. Therefore, we obtain $w \in L(\mathcal{S}(K))\phi$, which implies that $\overline{u} = \overline{w} \in \overline{L(\mathcal{S}(K))\phi} = \overline{K\varphi}$, concluding the proof.
\end{proof}

We begin by analyzing the effect of an automorphism of type $1$ at the level of Stallings automata. The following lemma is just a trivial observation.

\begin{lem}\label{elem:monoid:same}
Let $K \leq_{f.g.} F_A$. Then $M(K\alpha_{a}) = M(K)$ for all $a \in A$.
\end{lem}

Given a class of monoids $\mathcal{C}$ and $\varphi \in \mathrm{Aut}(F_A)$, we say that $\mathcal{C}$ is \textbf{stable under $\varphi$} if the condition \[ M(K) \in \mathcal{C} \Rightarrow M(K\varphi) \in \mathcal{C}\]
holds for all $K \leq_{f.g.} F_A$. 

The next lemma contains two other simple remarks.
\begin{lem}\label{compos:elem}
Let $\mathcal{C}$ be any class of monoids. 
\begin{enumerate}[label=(\roman*)]\setlength\itemsep{-0.1em}
\item If $\mathcal{C}$ is stable under all automorphisms of types $1$ and $2$, then $\mathcal{C}$ is stable under all automorphisms of $F_A$.
\item If $M(K\varphi) = M(K)$ holds for all automorphisms $\varphi$ of type 1 and 2, then it holds for all automorphisms of $F_A$.
\end{enumerate}
\end{lem}

Given $K \leq_{f.g.} F_A$, let $\mathcal{S}(K) = (Q,A,\delta,q_0)$ be the Stallings automaton of $K$. We proceed by analyzing the effect of an automorphism of type $2$ on $\St(K)$. Out of convenience, we paint the vertices in $Q$ \red{red}. 

Let $a, b \in A$ be two distinct letters, and write $\beta = \beta_{ab}$. When no confusion arises, we will always write $\beta$ instead of $\beta_{ab}$.
Consider the automaton $\mathcal{B}$ obtained from $\mathcal{S}(K)$ by replacing each edge 
$q$ \hspace{-2mm} $\xrightarrow{\; \, a \; \,}$ \hspace{-2mm} $q'$ (and its inverse)
with 
$q$ \hspace{-2mm} $\xrightarrow{\; \, a \; \,}$ \hspace{-2mm} $r$ \hspace{-2mm} $\xrightarrow{\; \, b \; \,}$ \hspace{-2mm} $q'$ (and their inverses), where $r \not \in Q$ is a new vertex which we paint \blue{blue}. This way, $\mathcal{B}$ becomes an involutive automaton. Observe that, by doing this, $\mathcal{B}$ is not necessarily deterministic, as we may encounter something like
\[ \xymatrix{\blue{\bullet} \ar[r]^-b & \red{\bullet} & \red{\bullet} \ar[l]_-b}.\]
In fact, ambiguity only arises in these cases, which appear if and only if we find 
\[ \xymatrix{\red{\bullet} \ar[r]^-a & \red{\bullet} & \red{\bullet} \ar[l]_-b}\]
in $\mathcal{S}(K)$. After we perform all the foldings in $\mathcal{B}$ and delete all the vertices with outdegree $1$ which are not the basepoint, we get an inverse automaton which is (isomorphic to) $\mathcal{S}(K\beta)$, in view of Proposition~\ref{prop:edge:replace}.

An important observation is that by doing only ``first order'' foldings, i.e., by replacing each ocurrence of 
\[ \xymatrix{\red{q} \ar[r]^-a & \blue{\bullet} \ar[r]^-b & \red{q'} & \red{q''} \ar[l]_-b }\]
in $\mathcal{B}$ with
\[ \xymatrix{\red{q} \ar[r]^-a & \red{q''} \ar[r]^-b & \red{q'}},\]
we already get a deterministic automaton. For a proof, see \cite{tese}.

Alternatively, consider the involutive automaton $\mathcal{S}(K)^\beta = (Q^\beta, A, \theta, q_0)$ obtained from $\mathcal{S}(K)$ by doing the following:
\begin{enumerate}\setlength\itemsep{-0.1em}
\item We keep all (positive) edges with label $x \in A \setminus\{a\}$.
\item If $q$ \hspace{-2mm} $\xrightarrow{\; \, a \; \,}$ \hspace{-2mm} $q'$ \hspace{-2mm} $\xleftarrow{\; \, b \; \,}$ \hspace{-2mm} $q''$ appears in $\mathcal{S}(K)$, we replace $q$ \hspace{-2mm} $\xrightarrow{\; \, a \; \,}$ \hspace{-2mm} $q'$ with $q$ \hspace{-2mm} $\xrightarrow{\; \, a \; \,}$ \hspace{-2mm} $q''$.
\item If $q$ \hspace{-2mm} $\xrightarrow{\; \, a \; \,}$ \hspace{-2mm} $q'$ appears in $\mathcal{S}(K)$, but it is not possible to read $b^{-1}$ from $q'$, then we replace that edge with $q$ \hspace{-2mm} $\xrightarrow{\; \, a \; \,}$ \hspace{-2mm} $r$ \hspace{-2mm} $\xrightarrow{\; \, b \; \,}$ \hspace{-2mm} $q'$, where $r \not \in Q$ is a new vertex.
\end{enumerate}

Note that $Q \subseteq Q^\beta$ and $\mathcal{S}(K)^\beta$ is isomorphic to the automaton we get after performing all ``first order'' foldings in $\mathcal{B}$, and hence it is an inverse automaton since no ``second order'' foldings are required. We keep considering that the vertices in $Q$ are painted red, and the ones in $Q^\beta \setminus Q$ are painted blue. 
It follows that $\mathcal{S}(K\beta) = (Q',A,\gamma,q_0)$ is the (inverse) automaton we obtain from $\mathcal{S}(K)^\beta$ by deleting all vertices with outdegree $1$ which are not the basepoint.

Observe that $\mathcal{S}(K\beta^{-1})$ can be obtained by a similar procedure. In this case, the notation we use is $\mathcal{S}(K)^{\beta^{-1}} = (Q^{\beta^{-1}},A,\theta',q_0)$.

\begin{exmp}
\normalfont Let $A = \{a,b,c\}$ and let $K \leq_{f.g.} F_A$ be given by
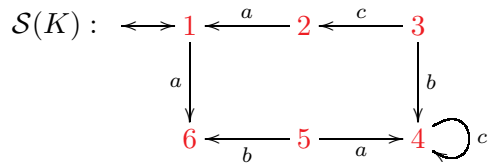
\[\St(K): \xygraph{  
!{<0cm,0cm>;<1cm,0cm>:<0cm,1cm>::}  
!{(-1,0)}*+{}="-1" 
!{(0,0)}*+{\red{1}}="1"    
!{(1.5,0)}*+{\red{2}}="2" 
!{(3,0)}*+{\red{3}}="3" 
!{(3,-1.5)}*+{\red{4}}="4" 
!{(1.5,-1.5)}*+{\red{5}}="5" 
!{(0,-1.5)}*+{\red{6}}="6"
"-1":"1" "1":"-1"  
"2":"1"_-{a}
"3":"2"_-{c}
"3":"4"^-{b}
"5":"4"_-{a}
"5":"6"^-{b}
"1":"6"_-{a}
"4":@(ru,rd)"4"^-{c}
} \]
The automaton $\B$ described above is obtained by replacing every edge $\red{i} \xrightarrow{\; a \;} \red{j}$ with a pair of edges $\red{i} \xrightarrow{\; a \;} \blue{k} \xrightarrow{\; b \;} \red{j}$, as we depict below.

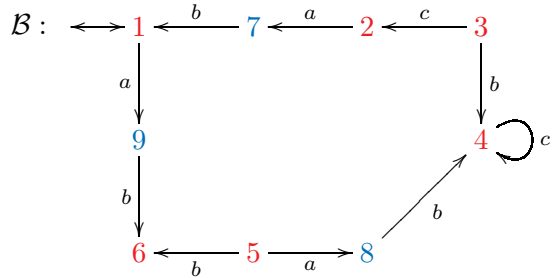
\[\B: \xygraph{  
!{<0cm,0cm>;<1cm,0cm>:<0cm,1cm>::}  
!{(-1,0)}*+{}="-1" 
!{(0,0)}*+{\red{1}}="1"  
!{(1.5,0)}*+{\blue{7}}="7" 
!{(3,0)}*+{\red{2}}="2" 
!{(4.5,0)}*+{\red{3}}="3" 
!{(4.5,-1.5)}*+{\red{4}}="4" 
!{(3,-3)}*+{\blue{8}}="8" 
!{(1.5,-3)}*+{\red{5}}="5" 
!{(0,-3)}*+{\red{6}}="6"
!{(0,-1.5)}*+{\blue{9}}="9"
"-1":"1" "1":"-1"  
"1":"9"_-{a}
"9":"6"_-{b}
"3":"2"_-{c}
"3":"4"^-{b}
"5":"8"_-{a}
"8":"4"_-{b}
"5":"6"^-{b}
"7":"1"_-{b}
"2":"7"_-{a}
"4":@(ru,rd)"4"^-{c}
} \]
After performing all the ``first order'' foldings in $\B$, we obtain an automaton isomorphic to $\St(K)^\beta$, which is also achieved by making the following substitutions in $\St(K)$: replace the edge $\red{2} \xrightarrow{\; a \;} \red{1}$ with $\red{2} \xrightarrow{\; a \;} \blue{7} \xrightarrow{\; b \;} \red{1}$; replace the edge $\red{5} \xrightarrow{\; a \;} \red{4}$ with $\red{5} \xrightarrow{\; a \;} \red{3}$; and replace the edge $\red{1} \xrightarrow{\; a \;} \red{6}$ with $\red{1} \xrightarrow{\; a \;} \red{5}$.

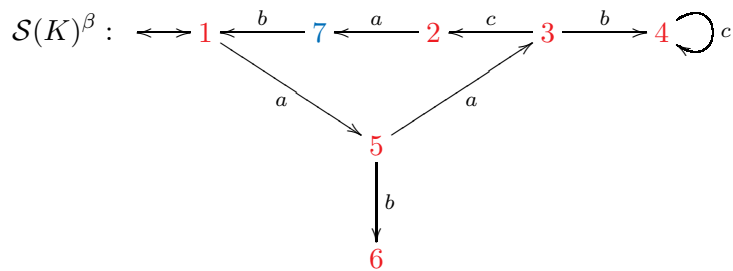
\[\St(K)^\beta: \xygraph{  
!{<0cm,0cm>;<1cm,0cm>:<0cm,1cm>::}  
!{(-1,0)}*+{}="-1" 
!{(0,0)}*+{\red{1}}="1"  
!{(1.5,0)}*+{\blue{7}}="7" 
!{(3,0)}*+{\red{2}}="2" 
!{(4.5,0)}*+{\red{3}}="3" 
!{(6,0)}*+{\red{4}}="4" 
!{(2.25,-1.5)}*+{\red{5}}="5" 
!{(2.25,-3)}*+{\red{6}}="6"
"-1":"1" "1":"-1"  
"1":"5"_-{a}
"3":"2"_-{c}
"3":"4"^-{b}
"5":"3"_-{a}
"5":"6"^-{b}
"7":"1"_-{b}
"2":"7"_-{a}
"4":@(ru,rd)"4"^-{c}
} \]
Finally, the Stallings automaton of $K\beta \leq_{f.g.} F_A$ is achieved by deleting vertex $6$ in $\St(K)^\beta$, because it has outdegree $1$ and it is not the basepoint.
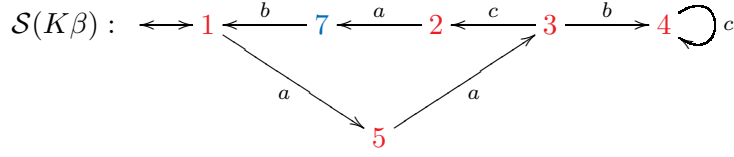
\[\St(K\beta): \xygraph{  
!{<0cm,0cm>;<1cm,0cm>:<0cm,1cm>::}  
!{(-1,0)}*+{}="-1" 
!{(0,0)}*+{\red{1}}="1"  
!{(1.5,0)}*+{\blue{7}}="7" 
!{(3,0)}*+{\red{2}}="2" 
!{(4.5,0)}*+{\red{3}}="3" 
!{(6,0)}*+{\red{4}}="4" 
!{(2.25,-1.5)}*+{\red{5}}="5" 
"-1":"1" "1":"-1"  
"1":"5"_-{a}
"3":"2"_-{c}
"3":"4"^-{b}
"5":"3"_-{a}
"7":"1"_-{b}
"2":"7"_-{a}
"4":@(ru,rd)"4"^-{c}
} \]
\end{exmp}

We now present some examples of pseudovarieties of monoids that \textbf{do not satisfy} the property of being stable under all automorphisms of a free group.

Let \textbf{Sl}, \textbf{R}, \textbf{L}, \textbf{Com} and \textbf{CR} be the pseudovarieties of semilattices, $\mathcal{R}$-trivial, $\mathcal{L}$-trivial, commutative and completely regular monoids, in that order. Let $A = \{a,b\}$, $\beta = \beta_{ab} \in \rr{Aut}(F_A)$ and $K = \langle a \rangle \leq F_A$ be the subgroup of $F_A$ generated by $a \in F_A$. Then the automata $\mathcal{S}(K) = (\{q_0\},\{a,b\},\delta,q_0)$ and $\mathcal{S}(K\beta) = (\{q_0,q_1\},\{a,b\},\gamma,q_0)$ are depicted by
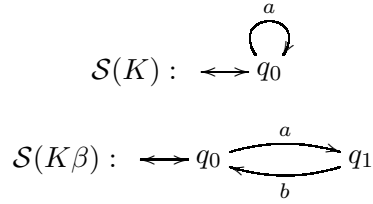
\[\mathcal{S}(K): \xygraph{  
!{<0cm,0cm>;<1cm,0cm>:<0cm,1cm>::}  
!{(-1,0)}*+{}="-1" 
!{(0,0)}*+{q_0}="1"    
"-1":"1" "1":"-1"  
"1":@(lu,ru)"1"^{a}
} \] 
\[\mathcal{S}(K\beta): \xygraph{  
!{<0cm,0cm>;<1cm,0cm>:<0cm,1cm>::}  
!{(-1,0)}*+{}="-1" 
!{(0,0)}*+{q_0}="1"    
!{(2,0)}*+{q_1}="2" 
"-1":"1" "1":"-1"  
"1":@/^/"2"^-{a} "2":@/^/"1"^-{b}
} \]  

We have the following:
\begin{itemize}\setlength\itemsep{-0.1em}
\item $M(K\beta) \not \in \boldsymbol{\mathrm{Sl}}$, for $\gamma_a$ is not an idempotent; however, $M(K) \in \boldsymbol{\mathrm{Sl}}$.
\item $M(K\beta) \not \in \boldsymbol{\mathrm{Com}}$, since $\gamma_a \gamma_b \neq \gamma_b \gamma_a$, but $M(K) \in \boldsymbol{\mathrm{Com}}$.
\item $M(K\beta) \not \in \boldsymbol{\mathrm{CR}} = \llbracket x^{\omega + 1} = x \rrbracket$, since $\gamma_a^\omega = 0 \neq \gamma_a $; in contrast, $M(K) \in \boldsymbol{\mathrm{CR}}$.
\end{itemize} 
These observations allow us to conclude that the pseudovarieties \textbf{Sl} (and, consequently, \textbf{R} and \textbf{L}, since every $\cl{R}$-trivial or $\cl{L}$-trivial inverse monoid is a semilattice), \textbf{Com} and \textbf{CR} are not stable under $\beta \in \mathrm{Aut}(F_A)$.

Moreover, despite the fact that the pseudovariety $\boldsymbol{\mathrm{A}} = \llbracket x^{\omega + 1} = x^\omega \rrbracket$ of aperiodic monoids is stable under all automorphisms of $F_A$, the pseudovariety $\boldsymbol{\mathrm{A}}_n = \llbracket x^{n+1} = x^n \rrbracket$ is not, for any $n \geq 1$. In fact, if $\xi \in \mathrm{Aut}(F_A)$ is the automorphism given by
\begin{align*}
\xi \colon F_A &\longrightarrow F_A \\
a &\longmapsto ab^n \\
b &\longmapsto b,
\end{align*}
then $M(K\xi) \not \in \boldsymbol{\mathrm{A}}_n$. To see why, observe that $\mathcal{S}(K\xi) = (\{q_0,\dots,q_n\},\{a,b\},\theta,q_0)$ can be depicted by
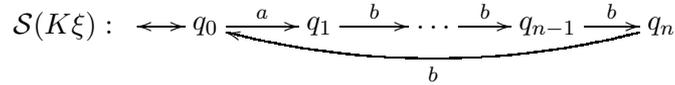
\[\mathcal{S}(K\xi): \xygraph{  
!{<0cm,0cm>;<1cm,0cm>:<0cm,1cm>::}  
!{(-1,0)}*+{}="-1" 
!{(0,0)}*+{q_0}="0"    
!{(1.5,0)}*+{q_1}="1" 
!{(3,0)}*+{\dots}="2"
!{(4.5,0)}*+{q_{n-1}}="3"
!{(6,0)}*+{q_n}="4"
"-1":"0" "0":"-1"  
"0":"1"^-{a} "1":"2"^-{b} "2":"3"^-{b} "3":"4"^-{b}
"4":@/^1pc/"0"^-{b} 
} \] 
So $\theta_b^{n+1} \neq \theta_b^n$, which shows that $M(K\xi) \not \in \boldsymbol{\mathrm{A}}_n$, even though $M(K) \in \boldsymbol{\mathrm{A}}_n$. 

We end this section with a lemma that we will use extensively. Let $\beta_0 \colon {\widetilde{A}}^* \longrightarrow {\widetilde{A}}^*$ be the unique free monoid homomorphism extending the function
\begin{align*}
\widetilde{A} &\longrightarrow {\widetilde{A}}^* \\
a &\longmapsto ab \\
a^{-1} & \longmapsto b^{-1}a^{-1} \\
x &\longmapsto x, \; \text{if } x \in \widetilde{A} \setminus \{a,a^{-1}\}.
\end{align*}
Similarly, let $\beta^{-1}_0 \colon {\widetilde{A}}^* \longrightarrow {\widetilde{A}}^*$ be the extension of
\begin{align*}
\widetilde{A} &\longrightarrow {\widetilde{A}}^* \\
a &\longmapsto ab^{-1} \\
a^{-1} & \longmapsto ba^{-1} \\
x &\longmapsto x, \; \text{if } x \in \widetilde{A} \setminus \{a,a^{-1}\}.
\end{align*}
We have the following result.

\begin{lem}\label{ab:lema}
Let $q,q' \in Q$ and $u \in {\widetilde{A}}^*$. Given $K \leq_{f.g.} F_A$, let $\St(K) = (Q,A,\delta,q_0)$ and $\St(K)^\beta = (Q^\beta,A,\theta,q_0)$.
\begin{enumerate}[label=\textup{(\alph*)}]\setlength\itemsep{-0.1em}
\item If $q \cdot u = q'$ in $\St(K)$, then $q \cdot u\beta_0 = q'$ in $\St(K)^\beta$. As a consequence, $q \cdot \overline{u\beta_0} = q'$ in $\St(K)^\beta$.
\item If $q \cdot u = q'$ in $\St(K)^\beta$, then $q \cdot \overline{u\beta^{-1}_0} = q'$ in $\St(K)$.
\end{enumerate}
\end{lem}

\begin{proof}
If $u=1$, both items are trivial, so we carry on with $u \in \widetilde{A}^+$.
\begin{enumerate}[label=\textup{(\alph*)}]
\item We can write $u = u_1u_2 \dots u_m$, where each $u_i$ is either a letter, the word $ab^{-1}$ or the word $ba^{-1}$, such that $u_iu_{i+1} \not \in \{ ab^{-1}, ba^{-1} \}$, for $1 \leq i \leq m-1$. Indeed, such factorization of $u$ is always possible, for $ab^{-1}$ and $ba^{-1}$ have no letters in common, and it is unique. We then have a path 
\[ q = \red{p_0} \xrightarrow{u_1} \red{p_1} \xrightarrow{u_2} \dots \xrightarrow{u_{m-1}} \red{p_{m-1}} \xrightarrow{u_m} \red{p_m} = q' \] in $\mathcal{S}(K)$, for some vertices $p_i \in Q$ and $0 \leq i \leq m$.

Given any $1 \leq i \leq m$, we analyse three possible cases. 
\begin{itemize}\setlength\itemsep{-0.1em}
\item If $u_i \not \in \left\{a^\varepsilon, {(ab^{-1})}^\varepsilon \mid \varepsilon = \pm 1 \right\}$, then \[ (p_{i-1},u_i)\delta = (p_{i-1},u_i)\theta = (p_{i-1},u_i\beta_0)\theta.\]
\item If $u_i = a^\varepsilon$ with $\varepsilon \in \{-1,1\}$, then \[(p_{i-1},u_i)\delta = (p_{i-1},a^\varepsilon)\delta = \left(p_{i-1},(ab)^\varepsilon\right)\theta = \left(p_{i-1},u_i \beta_0\right)\theta.\]
\item If $u_i = \left(ab^{-1}\right)^\varepsilon$ with $\varepsilon \in \{-1,1\}$, then \[(p_{i-1},u_i)\delta = \left(p_{i-1},(ab^{-1})^\varepsilon\right)\delta = (p_{i-1},a^\varepsilon)\theta = \left(p_{i-1},(abb^{-1})^\varepsilon\right)\theta = (p_{i-1},u_i \beta_0)\theta,\] since $\theta_{a}= \theta_{abb^{-1}}$, by construction of $\mathcal{S}(K)^\beta$.
\end{itemize}

To finish, one easily shows by induction on $m$ that $(q,u_1 \dots u_m)\delta = (q,(u_1 \dots u_m) \beta_0)\theta$, which concludes the proof.

\item Note that $(q,u)\delta = q' \Rightarrow (q,\overline{u})\delta = q'$ and $\overline{u\beta^{-1}_0} = \overline{\overline{u}\beta^{-1}_0}$, so we can assume that $u$ is reduced. The proof is similar to the one above. Indeed, we can write $u = u_1u_2 \dots u_m$, where each $u_i$ is either a letter, the word $ab$ or the word $b^{-1}a^{-1}$, such that $u_iu_{i+1} \not \in \{ ab, b^{-1}a^{-1} \}$, for $1 \leq i \leq m-1$. We get a path \[ q = \red{p_0} \xrightarrow{u_1} \red{p_1} \xrightarrow{u_2} \dots \xrightarrow{u_{m-1}} \red{p_{m-1}} \xrightarrow{u_m} \red{p_m} = q' \] in $\mathcal{S}(K)^\beta$, for some vertices $p_i \in Q$ and $0 \leq i \leq m$.

Let $\varepsilon \in \{-1,1\}$. Given any $1 \leq i \leq m$, again there are three possible cases.
\begin{itemize}\setlength\itemsep{-0.1em}
\item If $u_i \not \in \left\{a^\varepsilon, {(ab)}^\varepsilon\right\}$, then \[(p_{i-1},u_i)\theta = (p_{i-1},u_i)\delta = (p_{i-1},u_i\beta^{-1}_0)\delta = (p_{i-1},\overline{u_i\beta^{-1}_0})\delta.\]
\item If $u_i = a^\varepsilon$, then \[(p_{i-1},u_i)\theta = (p_{i-1},a^\varepsilon)\theta = \left(p_{i-1},(ab^{-1})^\varepsilon\right)\delta = \left(p_{i-1},u_i \beta^{-1}_0\right)\delta = \left(p_{i-1},\overline{u_i \beta^{-1}_0}\right)\theta.\]
\item If $u_i = \left(ab\right)^\varepsilon$, then \[(p_{i-1},u_i)\theta = \left(p_{i-1},\left(ab\right)^\varepsilon\right)\theta = (p_{i-1},a^\varepsilon)\delta = (p_{i-1},\overline{u_i \beta^{-1}_0})\delta,\] since $u_i\beta^{-1}_0 = \left(ab^{-1}b\right)^\varepsilon$.
\end{itemize}

Again, by induction on $m$, one shows that $(q,u_1 \dots u_m)\theta = (q,\overline{u_1 \beta^{-1}_0} \dots \overline{u_m \beta^{-1}_0})\delta$, which concludes the proof.
\end{enumerate}
\end{proof}

When dealing with the automorphism $\beta^{-1}$ instead of $\beta$, by switching the roles of $\beta_0$ and $\beta^{-1}_0$ we get the following analogue of Lemma \ref{ab:lema}.

\begin{lem}\label{ab^-1:lema}
Let $q,q' \in Q$ and $u \in {\widetilde{A}}^*$. Given $K \leq_{f.g.} F_A$, let $\St(K) = (Q,A,\delta,q_0)$ and $\St(K)^{\beta^{-1}} = (Q^{\beta^{-1}},A,\theta',q_0)$.
\begin{enumerate}[label=\textup{(\alph*)}]\setlength\itemsep{-0.1em}
\item If $q \cdot u= q'$ in $\St(K)$, then  $q \cdot u\beta^{-1}_0 = q'$ in $\St(K)^{\beta^{-1}}$. As a consequence, $q \cdot \overline{u\beta^{-1}_0} = q'$ in $\St(K)^{\beta^{-1}}$.
\item If $q \cdot u = q'$ in $\St(K)^{\beta^{-1}}$, then $q \cdot \overline{u\beta_0} = q'$ in $\St(K)$.
\end{enumerate}
\end{lem}

\section{Transition groups}\label{sec:transgroup}

In this small section, we focus on finitely generated subgroups $K \leq_{f.g.} F_A$ of a free group with the property that $M(K)$ is a group. We begin with a simple lemma.

\begin{lem}\label{lema:complete}
Let $\mathcal{A}$ be an inverse automaton over the alphabet $A$ which is not complete. Then there exists a word $w \in \widetilde{A}^*$ which cannot be read from any state in $\A$. 
\end{lem}

\begin{proof}
Let $\A$ be an inverse and incomplete automaton, with set of states $Q = \{q_1, \dots, q_n\}$ for some $n \geq 1$, and transition function $\delta$. We can assume that there exists a letter $a \in A$ such that $q_1 \not \in \dom \delta_a$. Since $\A$ is inverse, for each $1 \leq i \leq n$ there exists $u_i \in \widetilde{A}^*$ such that $q_i \cdot u_i = q_1$. Consider the word \[w = (u_1 a)(u_1 a)^{-1}(u_2 a)(u_2 a)^{-1}\dots (u_n a)(u_n a)^{-1}.\]
Since every factor of the type $(u_j a)(u_j a)^{-1}$ labels a loop at every state in $\dom \delta_{(u_j a)(u_j a)^{-1}}$, and $q_i \not \in \dom \delta_{u_ia}$ for every $1 \leq i \leq n$, it follows that $w$ is a word which cannot be read from any state in $\A$. 
\end{proof}

The next result is well known and a proof can be found in \cite[Theorem~5.1]{margolis}.

\begin{thm}\label{thm:completegroup}
Consider $1 \neq K \leq_{f.g.} F_A$ with Stallings automaton $\St(K) = (Q,A,\delta,q_0)$. The following conditions are equivalent:
\begin{enumerate}[label=(\roman*)]\setlength\itemsep{-0.1em}
\item $K$ has finite index in $F_A$, i.e., $[F_A:K] < \infty$.
\item $\St(K)$ is complete.
\item $M(K)$ is a group.
\end{enumerate}
\end{thm}

Now we prove that every class of groups is stable under all automorphisms of a free group.

\begin{thm}\label{groups:ab}
Let $K \leq_{f.g.} F_A$ and let $\varphi$ be an automorphism of $F_A$. If $M(K)$ is a group, then $M(K\varphi) = M(K)$.
\end{thm}

\begin{proof}
The case $K=1$ is clear, so we focus on the case $K \neq 1$. Suppose that $M(K)$ is a group, let $\mathcal{S}(K) = (Q,A,\delta,q_0)$ be the Stallings automaton of $K$ and $\mathcal{S}(K\beta) = (Q',A,\gamma,q_0)$ be the Stallings automaton of $K\beta$. Since $\mathcal{S}(K)$ is complete by Theorem~\ref{thm:completegroup}, $\mathcal{S}(K)^\beta$ is the automaton obtained from $\mathcal{S}(K)$ by simply replacing every occurrence of $q$ \hspace{-2mm} $\xrightarrow{\; \, a \; \,}$ \hspace{-2mm} $q'$ \hspace{-2mm} $\xleftarrow{\; \, b \; \,}$ \hspace{-2mm} $q''$ with $q$ \hspace{-2mm} $\xrightarrow{\; \, a \; \,}$ \hspace{-2mm} $q''$ \hspace{-2mm} $\xrightarrow{\; \, b \; \,}$ \hspace{-2mm} $q'$. Moreover, since $\mathcal{S}(K)^\beta$ is also complete, it has no vertices with outdegree $1$, whence $\mathcal{S}(K\beta) = \mathcal{S}(K)^\beta$.

Now, we know that $M(K\beta)$ is generated, as a group, by $\gamma_a$ and all the $\gamma_x$ with $x \in A \setminus \{a\}$. But $M(K\beta)$ is also generated by $\gamma_{ab}$ and $\gamma_x$ ($x \in A \setminus \{a\}$), since we can write $\gamma_{a} = \gamma_{ab} {\gamma_b}^{-1}$. Moreover, note that $\gamma_{ab} = \delta_a$ and $\gamma_x = \delta_x$ for all $x \in A \setminus \{a\}$. Since $\delta_a$, $\delta_x$ ($x \in A \setminus \{a\}$) generate $M(K)$ as a group, we get $M(K) = M(K\beta)$.
By a similar argument, we get $M(K) = M(K\beta^{-1})$. The conclusion now follows from Lemma~\ref{elem:monoid:same} and Lemma~\ref{compos:elem}.
\end{proof}

When the transition monoid of $\St(K)$ is a group, we can actually characterize $K$ in terms of the pseudovarieties $M(K)$ belongs to.

\begin{thm}\label{thm:M(K)group}
Let $K \leq_{f.i.} F_A$ with $\St(K) = (Q,A,\delta,q_0)$. Consider the pseudovariety of groups $\boldsymbol{\mathrm{V}}$ ultimately defined by the sequence of identities $(u_n = 1)_n$ over the set of variables $X$. The following conditions are equivalent:
\begin{enumerate}[label=(\roman*)]\setlength\itemsep{-0.1em}
\item $M(K) \in \boldsymbol{\mathrm{V}}$.
\item There exists $p\geq 1$ such that, for every group homomorphism $\psi \colon F_X \rightarrow F_A$, the condition $u_n\psi \in K$ holds for all $n \geq p$.
\end{enumerate}
\end{thm}

\begin{proof}
(i) $\Rightarrow$ (ii): By assumption, $K$ has finite index in $F_A$, so the Stallings automaton $\mathcal{S}(K)$ is complete. This implies that if $v,w \in \widetilde{A}^*$ are such that $\overline{v} = \overline{w}$, then for every $q \in Q$ we have
\[q \cdot v = q \cdot \overline{v} = q \cdot \overline{w} = q \cdot w. \]
Hence, the monoid homomorphism
\begin{align*}
\Delta \colon {\widetilde{A}}^* &\longrightarrow M(K) \\
u &\longmapsto \delta_u
\end{align*}
induces a group homomorphism 
\begin{align*}
\widetilde{\Delta} \colon F_A &\longrightarrow M(K) \\
g &\longmapsto \delta_{\overline{g}}.
\end{align*}
By assumption, there exists $p \geq 1$ such that $M(K) \models u_n = 1$ for all $n \geq p$. This means that given any group homomorphism $\varphi \colon F_X \rightarrow M(K)$, we have $u_n \varphi = 1$ for all $n \geq p$. 
Fix such a $p \geq 1$ and let $\psi \colon F_X \rightarrow F_A$ be a homomorphism. Since $\psi \widetilde{\Delta} \colon F_X \rightarrow M(K)$ is a homomorphism, it follows that $u_n \psi \widetilde{\Delta} = 1$ for all $n \geq p$. This implies that $\delta_{\overline{u_n\psi}} = 1$ for all $n \geq p$, and thus $q \cdot \overline{u_n\psi} = q$ for all $q \in Q$ and $n \geq p$. In particular, we have $q_0 \cdot \overline{u_n\psi}=q_0$ for all $n \geq p$, which shows that $u_n \psi \in K$ for all $n \geq p$, as we wanted.

(ii) $\Rightarrow$ (i): Let $\phi \colon \widetilde{A}^* \rightarrow F_A$ be the quotient morphism. Take $p \geq 1$ satisfying (ii) and $n \geq p$. There exist $x_{1}, \dots, x_{m} \in X$ and $k_{1}, \dots, k_{m} \in \mathbb{Z}$ such that $u_n = x_{1}^{k_1} \cdots x_{m}^{k_m} \in F_X$. Our goal is to prove that, given any words $w_1, \dots, w_m \in \widetilde{A}^*$, we have $\delta_{w_1}^{k_1} \cdots \delta_{w_m}^{k_m} = 1$, which is the same as proving that $q \cdot w_1^{k_1} \cdots w_m^{k_m} = q$ for all $q \in Q$. Note that the condition $K \leq_{f.i.} F_A$ entails that $\St(K)$ is complete.

Write $w = w_1^{k_1}\cdots w_m^{k_m}$, take any $q \in Q$, and let $v \in \widetilde{A}^*$ be such that $q_0 \cdot v = q$.
Since $F_X$ is a free group, the function 
$X \rightarrow F_A$, $x_i \mapsto (vw_iv^{-1})\phi$
induces a homomorphism $\psi \colon F_X \rightarrow F_A$ satisfying $u_n \psi = (vwv^{-1})\phi$. Then $(vwv^{-1})\phi \in K$ and thus $q_0 \cdot vwv^{-1} = q_0$, by construction of $\mathcal{S}(K)$ and taking into account that it is a complete automaton. We conclude that $q \cdot w = q$, as we wanted.
\end{proof}

\begin{exmp}
\normalfont The pseudovariety of abelian groups is given by $\br{Ab} = \llbracket xyx^{-1}y^{-1} = 1 \rrbracket$, and the pseudovariety of $p$-groups is given by $\br{G}_p = \llbracket x^{p^\omega} = 1 \rrbracket$. Hence, the above theorem allows us to conclude that $M(K)$ is an abelian group if and only if $K$ has finite index and contains all elements of the type $xyx^{-1}y^{-1}$ $(x,y \in F_A)$; and it is a $p$-group if and only if it has finite index and there exists $m \geq 1$ such that $x^{p^{m}} \in K$ for all $x \in F_A$.
\end{exmp}

\section{Pseudovarieties of type \texorpdfstring{$\overline{\br{H}}$}{H}}\label{sec:main}

In this section, we prove that if $\br{H}$ is a pseudovariety of finite groups, then the pseudovariety of monoids $\overline{\br{H}}$ (defined at the end of Section~\ref{subsec:psv}) is stable under all automorphisms of $F_A$. Then we characterize some finitely generated subgroups of $F_A$ in terms of these pseudovarieties.

We begin with a simple lemma regarding the group $\mathcal{H}$-classes of an inverse monoid, which follows from \cite[Lemma II.1.7]{petrich}.

\begin{lem}\label{lema:group:Hclass}
Let $M$ be an inverse monoid and denote by $H_x$ the $\mathcal{H}$-class of $x \in M$. Then $H_x$ is a group if and only if $xx^{-1} = x^{-1}x$. Moreover, $y \in H_x$ if and only if $yy^{-1} = xx^{-1}$ and $y^{-1}y = x^{-1}x$.
\end{lem}

Observe that if the elements of a monoid $M$ are injective partial transformations, then $f,g \in M$ are in the same group $\mathcal{H}$-class if and only if $\dom f = \im f = \dom g = \im g$. 

Recall that $\St(K) = (Q,A,\delta,q_0)$, $\St(K)^\beta = (Q^\beta,A,\theta,q_0)$ and let $M(K)^\beta$ denote the transition monoid of $\mathcal{S}(K)^\beta$. We will now show that every group $\mathcal{H}$-class of $M(K)^\beta$ is the homomorphic image of a subgroup of some group $\mathcal{H}$-class of $M(K)$. We say that a monoid $M$ divides a monoid $N$ if $M$ is a homomorphic image of a submonoid of $N$.

\begin{lem}\label{lemma:Hclass:theta}
Every group $\mathcal{H}$-class of $M(K)^\beta$ divides some group $\mathcal{H}$-class of $M(K)$.
\end{lem}

\begin{proof}
We invite the reader to follow the arguments on Example~\ref{ex:h}. Let $H$ be a group $\mathcal{H}$-class of $M(K)^\beta$. We will analyse two cases and make use of Lemma \ref{ab:lema} several times without explicitly mentioning it. 

\begin{itemize}
\item \textbf{Case 1:} There is some $u \in \widetilde{A}^+$ such that $\mathrm{dom}\, \theta_{uu^{-1}} \subseteq Q$ and $H = H_{\theta_{uu^{-1}}}$.

Let $X = \mathrm{dom} \, \theta_{uu^{-1}} = \mathrm{dom} \, \theta_{u^{-1}u}$ and define \[ X' = \bigcap\limits_{\substack{v \in \widetilde{A}^* \\ \theta_v \in H}} \mathrm{dom}\, \delta_{\overline{v \beta^{-1}_0}}. \]
Note that $X \subseteq X'$ because, given $q \in X$, we have $q \in \mathrm{dom}\, \theta_v$ for every $v \in \widetilde{A}^*$ satisfying $\theta_v \in H$, and thus $q \in \mathrm{dom}\, \delta_{\overline{v \beta^{-1}_0}}$, which implies $q \in X'$. Furthermore, since $\mathcal{S}(K)$ is a finite automaton, there exists some $\ell \geq 0$ such that, for every $w \in \widetilde{A}^*$, we can find a word $\widetilde{w}$ satisfying $|\widetilde{w}| \leq \ell$ and $\delta_w = \delta_{\tilde{w}}$. Hence, there is a finite amount of words $v_1, v_2, \dots, v_m \in \widetilde{A}^*$ ($m \geq 0$) such that $\theta_{v_1}, \theta_{v_2}, \dots, \theta_{v_m} \in H$ and \[X' = \bigcap_{i=1}^m \mathrm{dom}\, \delta_{\overline{v_i \beta^{-1}_0}}.\] Consider the element $e = \overline{v_1 \beta^{-1}_0} \, {\overline{v_1 \beta^{-1}_0}}^{-1} \dots \overline{v_m \beta^{-1}_0} \, {\overline{v_m \beta^{-1}_0}}^{-1}$. Then $\delta_e$ is an idempotent and $\mathrm{dom}\, \delta_e = X'$. Now, consider the subgroup $H' \leq H_{\delta_e}$ given by 
\[ \delta_w \in H' \Leftrightarrow X\delta_w \subseteq X. \]
In fact, since $\delta_w \in H'$ is injective, we have $X\delta_w \subseteq X \Leftrightarrow X\delta_w = X$; and if $\delta_w$ satisfies this condition, then also $(X' \setminus X)\delta_w = X' \setminus X$, since $\mathrm{dom}\, \delta_w = \mathrm{dom}\, \delta_{w^{-1}}$. We claim that the function $\Phi$ given by
\begin{align*}
\Phi \colon H' &\longmapsto H \\
\delta_w &\longrightarrow \theta_{uu^{-1}}\theta_{w\beta_0}
\end{align*}
is a surjective group homomorphism. 

First we show that $\Phi$ is well defined. Given $\delta_w \in H'$, we must prove that 
\begin{equation}\label{eq:1}
\mathrm{dom }\, (\theta_{uu^{-1}}\theta_{w\beta_0}) = X
\end{equation}
\vspace{-8mm}
\begin{equation}\label{eq:2}
\mathrm{dom }\, (\theta_{w^{-1}\beta_0}\theta_{uu^{-1}}) = X.
\end{equation}
We start with (\ref{eq:1}). It is clear that the left hand side of the equality is contained in $X$. Regarding the reverse inclusion, given $q \in X$, we have $q \in X' = \mathrm{dom}\, \delta_{w}$, which implies $q \in \mathrm{dom} \, \theta_{w\beta_0}$ and $q \in \mathrm{dom }\, (\theta_{uu^{-1}}\theta_{w\beta_0})$. Now we focus on (\ref{eq:2}). If $q \in X$ then $q \in \mathrm{dom} \, \delta_{w^{-1}}$, and thus $q \in \mathrm{dom} \, \theta_{w^{-1}\beta_0}$. Consequently, \[(q, w^{-1}\beta_0)\theta = (q,w^{-1})\delta \in X = \mathrm{dom}\, \theta_{uu^{-1}}\] because $X\delta_{w^{-1}} = X$. As for the other inclusion, if we take $q \in \mathrm{dom }\, (\theta_{w^{-1}\beta_0}\theta_{uu^{-1}})$, then \[q' = (q,w^{-1}\beta_0)\theta \in \mathrm{dom}\, \theta_{uu^{-1}} = X \subseteq \mathrm{dom} \, \delta_w.\] Therefore, we can read $w$ from $q'$ in $\mathcal{S}(K)$, and if $(q',w)\delta = q''$, then $(q',w\beta_0)\theta = q''$ by Lemma \ref{ab:lema}. Since the automaton $\mathcal{S}(K)^\beta$ is inverse, the equality $q' = (q,w^{-1}\beta_0)\theta$ implies $q = (q',w\beta_0)\theta$. It follows that $(q',w)\delta = q \in X$, as $X\delta_w = X$. This completes the proof of (\ref{eq:2}). 

In order for $\Phi$ to be indubitably well defined, we must also confirm that if $\delta_w = \delta_z \in H'$, then $\theta_{uu^{-1}}\theta_{w\beta_0} = \theta_{uu^{-1}}\theta_{z\beta_0}$. But we already know that \[\mathrm{dom} \, (\theta_{uu^{-1}}\theta_{w\beta_0}) = \mathrm{dom} \, (\theta_{uu^{-1}}\theta_{z\beta_0}) = X;\] and taking $q \in X \subseteq \mathrm{dom} \, \delta_w = \mathrm{dom} \, \delta_z$, we obtain
\begin{align*}
(q, uu^{-1}(w\beta_0))\theta = (q,w\beta_0)\theta = (q,w)\delta = (q,z)\delta = (q,z\beta_0)\theta = (q, uu^{-1}(z\beta_0))\theta,
\end{align*}
which establishes the claim.

Next, we show that $\Phi$ is a homomorphism, which amounts to proving that 
\[\theta_{uu^{-1}}\theta_{w\beta_0}\theta_{uu^{-1}}\theta_{z\beta_0} = \theta_{uu^{-1}}\theta_{w\beta_0}\theta_{z\beta_0}\]
for all $\delta_w, \delta_z \in H'$. Again, it suffices to check that these transformations have the same domain. Indeed, if that is true, then for any vertex $q$ in their domain, we will obtain
\begin{align*}
(q, uu^{-1}(w\beta_0)uu^{-1}(z\beta_0))\theta &= (((q,uu^{-1}(w\beta_0))\theta,uu^{-1})\theta,z\beta_0)\theta\\
&= ((q,uu^{-1}(w\beta_0))\theta,z\beta_0)\theta \\
&= (q, uu^{-1}(w\beta_0)(z\beta_0))\theta.
\end{align*}
It is clear that \[\mathrm{dom}\, (\theta_{uu^{-1}}\theta_{w\beta_0}\theta_{uu^{-1}}\theta_{z\beta_0}) \subseteq \mathrm{dom}\, (\theta_{uu^{-1}}\theta_{w\beta_0}\theta_{z\beta_0}).\] Regarding the other inclusion, let $q \in \mathrm{dom}\, (\theta_{uu^{-1}}\theta_{w\beta_0}\theta_{z\beta_0})$. Then $q \in X \subseteq \mathrm{dom}\,\delta_w$ and we can take $q' = (q,w\beta_0)\theta = (q,w)\delta \in X$, as $q \in X$ and $X\delta_w = X$. Since $q' \in \mathrm{dom}\, \theta_{z\beta_0}$, we conclude that $q' \in \mathrm{dom} \, (\theta_{uu^{-1}}\theta_{z\beta_0})$, whence $\Phi$ is indeed a group homomorphism.

Finally, we prove that $\Phi$ is surjective by showing that \[\theta_v = \theta_{uu^{-1}}\theta_{(e\overline{v \beta^{-1}_0}e)\beta_0} = (\delta_{e\overline{v \beta^{-1}_0}e})\Phi\] for every $\theta_v \in H$. Observe that $\delta_{e\overline{v \beta^{-1}_0}e} \in H_{\delta_e}$ because $\theta_{v^{\pm 1}} \in H$, $\overline{v\beta^{-1}_0}^{-1} = \overline{(v\beta^{-1}_0)^{-1}} = \overline{v^{-1}\beta^{-1}_0}$ and $\mathrm{dom }\, \delta_{e\overline{v^{\pm 1}\beta^{-1}_0}e} = X'$ by definition of $X'$ and $\delta_e$. Moreover, $\delta_{e\overline{v \beta^{-1}_0}e} \in H'$ since $q \in X$ implies \[(q, e\overline{v \beta^{-1}_0}e)\delta = (q, \overline{v \beta^{-1}_0})\delta = (q, v)\theta \in X. \] Hence, it remains to verify that $\theta_v = \theta_{uu^{-1}}\theta_{(e\overline{v \beta^{-1}_0}e)\beta_0}$. If these partial transformations have the same domain then, for any vertex $q$ in their domain, 
\begin{align*}
(q,uu^{-1}(e\overline{v\beta^{-1}_0}e)\beta_0)\theta &= (q, (e\beta_0) \overline{v\beta^{-1}_0} \beta_0 (e\beta_0))\theta \\
&= (q, \overline{v\beta^{-1}_0} \beta_0)\theta \\
&= (q, v)\theta.
\end{align*} 
Hence, it suffices to check that $\mathrm{dom}\, \theta_v = \mathrm{dom} \, (\theta_{uu^{-1}}\theta_{(e\overline{v \beta^{-1}_0}e)\beta_0})$. On the one hand, the equalities $\mathrm{dom}\, \theta_v = X = \mathrm{dom}\, \theta_{uu^{-1}}$ entail $\mathrm{dom}\, \theta_v \supseteq \mathrm{dom} \, (\theta_{uu^{-1}}\theta_{(e\overline{v \beta^{-1}_0}e)\beta_0})$. On the other hand, if $q \in \mathrm{dom}\, \theta_v = \mathrm{dom}\, \theta_{uu^{-1}}$, then $q \in \mathrm{dom}\, \delta_e$, which implies that $q \in \mathrm{dom}\, \theta_{e\beta_0}$. Also, $q \in \mathrm{dom}\, \delta_{\overline{v\beta^{-1}_0}}$ yields $q \in \mathrm{dom} \, \theta_{\overline{v\beta^{-1}_0}\beta_0}$. We conclude the proof by observing that $(q,\overline{v\beta^{-1}_0}\beta_0)\theta = (q,v)\theta \in X \subseteq \mathrm{dom} \, \delta_e \subseteq \mathrm{dom} \, \theta_{e\beta_0}$. 

\item \textbf{Case 2:} There is some $u \in \widetilde{A}^+$ such that $\mathrm{dom}\, \theta_{uu^{-1}} \not\subseteq Q$ and $H = H_{\theta_{uu^{-1}}}$.

We begin by recalling an important fact. Given a blue vertex $q \in Q^\beta \setminus Q$ in $\mathcal{S}(K)^\beta$ and a word $w \in \widetilde{A}^*$, if $q \in \mathrm{dom}\, \theta_w$, then the first letter of $w$ must be $a^{-1}$ or $b$; and if $q \in \mathrm{im}\, \theta_w$, then the last letter of $w$ must be $a$ or $b^{-1}$.

In this case, if $v \in \widetilde{A}^*$ is such that $\theta_v \in H$, then the equalities $\theta_a = \theta_{abb^{-1}}$ and $\theta_{a^{-1}} = \theta_{bb^{-1}a^{-1}}$ allow us to write $\theta_v = \theta_{b\tilde{v}b^{-1}}$ for some $\widetilde{v} \in \widetilde{A}^*$; in particular, $\theta_u = \theta_{b\tilde{u}b^{-1}}$. Now, since 
\[\theta_{b\tilde{u}b^{-1}b\tilde{u}^{-1}b^{-1}} \mathrel{\mathcal R} \theta_{(b\tilde{u}b^{-1}b\tilde{u}^{-1}b^{-1})b} \mathrel{\mathcal L} \theta_{b^{-1}(b\tilde{u}b^{-1}b\tilde{u}^{-1}b^{-1})b}, \]
we get $\theta_{uu^{-1}} \mathrel{\mathcal D} \theta_{(b^{-1}b\tilde{u}b^{-1})(b^{-1}b\tilde{u}b^{-1})^{-1}}$. In view of \cite[Proposition 2.3.6]{howie}, we conclude that  
\[H \cong H_{\theta_{b^{-1} b \tilde{u} b b^{-1} \tilde{u}^{-1} b^{-1}b}}.\] Therefore, this case can be reduced to the previous one.
\end{itemize}

Denoting by $\theta_1$ the transformation induced by the empty word $1$, observe that above we have covered all possibilites, since $H_{\theta_1}$ is either trivial (and the claim automatically follows), or there exists $x \in A$ satisfying $\theta_1 = \theta_{xx^{-1}}$ (and we fall into the previous cases).
\end{proof}

\begin{exmp}\label{ex:h}
\normalfont Let $A = \{a,b,c,d\}$ and let $K \leq_{f.g.} F_A$ be given by the Stallings automaton $\St(K) = (Q,A,\delta,1)$ depicted by
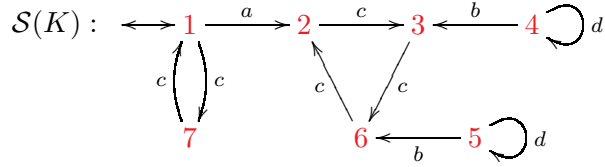
\[\St(K): \xygraph{  
!{<0cm,0cm>;<1cm,0cm>:<0cm,1cm>::}  
!{(-1,0)}*+{}="-1" 
!{(0,0)}*+{\red{1}}="1"  
!{(1.5,0)}*+{\red{2}}="2" 
!{(3,0)}*+{\red{3}}="3" 
!{(4.5,0)}*+{\red{4}}="4"
!{(3.75,-1.5)}*+{\red{5}}="5" 
!{(2.25,-1.5)}*+{\red{6}}="6"
!{(0,-1.5)}*+{\red{7}}="7"
"-1":"1" "1":"-1"  
"1":"2"^-{a} "2":"3"^-{c} "4":"3"_-{b} 
"5":"6"^-{b} "6":"2"^-{c} "3":"6"^-{c} 
"7":@/^/"1"^-{c} 
"1":@/^/"7"^-{c} 
"4":@(ru,rd)"4"^-{d}
"5":@(ru,rd)"5"^-{d}
} \]
Then $\St(K\beta) = \St(K)^\beta = (Q',A,\gamma,1)$ can be portrayed by
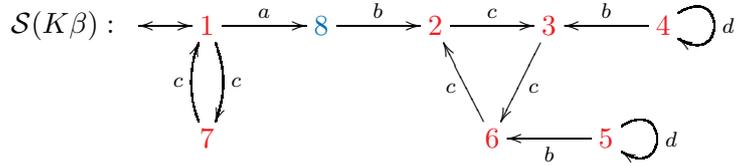
\[\St(K\beta): \xygraph{  
!{<0cm,0cm>;<1cm,0cm>:<0cm,1cm>::}  
!{(-1,0)}*+{}="-1" 
!{(0,0)}*+{\red{1}}="1"  
!{(1.5,0)}*+{\blue{8}}="8" 
!{(3,0)}*+{\red{2}}="2" 
!{(4.5,0)}*+{\red{3}}="3" 
!{(6,0)}*+{\red{4}}="4"
!{(5.25,-1.5)}*+{\red{5}}="5" 
!{(3.75,-1.5)}*+{\red{6}}="6"
!{(0,-1.5)}*+{\red{7}}="7"
"-1":"1" "1":"-1"  
"1":"8"^-{a} "8":"2"^-{b} "2":"3"^-{c} "4":"3"_-{b} 
"5":"6"^-{b} "6":"2"^-{c} "3":"6"^-{c} 
"7":@/^/"1"^-{c} 
"1":@/^/"7"^-{c} 
"4":@(ru,rd)"4"^-{d}
"5":@(ru,rd)"5"^-{d}
} \]
Let $u = b^{-1}bc$. Using the notation in the proof of Lemma~\ref{lemma:Hclass:theta}, consider the group $\cl{H}$-class
\[H = H_{\gamma_{uu^{-1}}} = \{ \gamma_u, \gamma_{u^2}, \gamma_{u^3}\} \cong \mathbb{Z}_3.\]
of $M(K\beta)$. Then we have 
\begin{align*}
X &= \dom \gamma_{uu^{-1}} = \{2,3,6\} \\
X' &= \dom \delta_{\overline{u \beta^{-1}_0}} \cap \dom \delta_{\overline{u^2\beta^{-1}_0}} \cap \dom \delta_{\overline{u^3 \beta^{-1}_0}} = \dom \delta_c \cap \dom \delta_{c^2} \cap \dom \delta_{c^3} = \{1,2,3,6,7\}.
\end{align*}
Letting $e = cc^{-1}$, it follows that
\begin{align*}
H_{\delta_e} &= H_{\delta_c} = \{\delta_c, \delta_{c^2}, \dots, \delta_{c^6}\} = H' \cong \mathbb{Z}_6
\end{align*}
Since $H \cong \mathbb{Z}_3$ is a homomorphic image of $H' \cong \mathbb{Z}_6$, one monoid divides the other.
\end{exmp}

We need one last lemma before proving the main theorem of this section. It can be proven in a similar way, so we only make a sketch here. For details, see \cite{tese}.

\begin{lem}\label{lem:reductions}
Let $\mathcal{A}' = (P',A,\lambda',q_0)$ be an inverse automaton with a basepoint and let $\mathcal{A} = (P,A,\lambda,q_0)$ be a subautomaton obtained by deleting a vertex of outdegree $1$ which is not the basepoint. Then every group $\mathcal{H}$-class of $M(\mathcal{A})$ divides some group $\mathcal{H}$-class of $M(\mathcal{A}')$.
\end{lem}

\begin{proof}
Let $P' \setminus P = \{p'\}$ and $p \xrightarrow{\; x \;} p'$, $p' \xrightarrow{\; x^{-1} \;} p$ ($p \in P$, $x \in \widetilde{A}$) be the extra edges in $\mathcal{A}'$. Let also $u \in \widetilde{A}^+$ be such that $H = H_{\lambda_{uu^{-1}}}$ is a group $\mathcal{H}$-class of $\mathcal{A}$. We define the sets $Y = \mathrm{dom} \, \lambda_{uu^{-1}} = \mathrm{dom} \, \lambda_{u^{-1}u}$ and \[ Y' = \bigcap\limits_{\substack{v \in \widetilde{A}^* \\ \lambda_v \in H}} \mathrm{dom}\, \lambda'_v. \]
It is clear that $Y \subseteq Y'$, and we can choose a finite amount of words $v_1, v_2, \dots, v_m \in \widetilde{A}^*$ ($m \geq 0$) such that $\lambda_{v_1}, \lambda_{v_2}, \dots, \lambda_{v_m} \in H$ and $Y' = \bigcap_{i=1}^m \mathrm{dom}\, \lambda'_{v_i}$. 
Consider the word \[ e = v_1v_1^{-1} v_2 v_2^{-1} \dots v_m v_m^{-1} \in \widetilde{A}^*.\] Then $\lambda'_e$ is an idempotent and $\mathrm{dom}\, \lambda'_e = Y'$. 
As in the previous lemma, we split the proof in two cases and argue in a similar manner.

\begin{itemize}
\item \textbf{Case 1:} $p' \not \in Y'$.

Let $H' \leq H_{\lambda'_e}$ be the subgroup defined by 
\[ \lambda'_w \in H' \Leftrightarrow Y\lambda'_w = Y. \]
Then it suffices to check that
\begin{align*}
\Psi \colon H' &\longmapsto H \\
\lambda'_w &\longrightarrow \lambda_{uu^{-1}}\lambda_{\overline{w}}
\end{align*}
is a surjective group homomorphism. 

\item \textbf{Case 2:} We have $p' \in Y'$.

In this case, if $v \in \widetilde{A}^+$ is such that $\lambda_v \in H$ then, since $p' \in \mathrm{dom}\, \lambda'_v \cap \mathrm{dom}\, \lambda'_{v^{-1}}$ has outdegree $1$, we must have $v = x^{-1} \widetilde{v} x$ for some $\widetilde{v} \in \widetilde{A}^*$. In particular, $u = x^{-1} \widetilde{u} x$ for some $\widetilde{u} \in \widetilde{A}^*$. Similarly to what we obtained in Case 2 of Lemma~\ref{lemma:Hclass:theta}, it follows that
\[\lambda_{x^{-1}\tilde{u}xx^{-1}\tilde{u}^{-1}x} \mathrel{\mathcal R} \lambda_{(x^{-1}\tilde{u}xx^{-1}\tilde{u}^{-1}x)x^{-1}} \mathrel{\mathcal L} \lambda_{x(x^{-1}\tilde{u}xx^{-1}\tilde{u}^{-1}x)x^{-1}}, \]
which implies $\lambda_{uu^{-1}} \mathrel{\mathcal D} \lambda_{(xx^{-1}\tilde{u}x)(xx^{-1}\tilde{u}x)^{-1}}$ and
\[H \cong H_{\lambda_{xx^{-1} \tilde{u} x^{-1} x \tilde{u}^{-1} x x^{-1}}}.\] This case can now be reduced to the previous one, as we cannot read a word starting with $x$ from $p'$.
\end{itemize}
\end{proof}

As we have remarked before, the automaton $\mathcal{S}(K\beta) = (Q',A,\gamma,q_0)$ can be obtained from $\mathcal{S}(K)^\beta = (Q^\beta,A,\theta,q_0)$ by successively deleting all vertices in $Q^\beta \setminus \{q_0\}$ with outdegree $1$. Hence, in view of Lemma~\ref{lem:reductions}, it follows that every group $\mathcal{H}$-class of $M(K\beta)$ divides some group $\mathcal{H}$-class of $M(K)^\beta$. Applying then Lemma~\ref{lemma:Hclass:theta}, we conclude that every group $\mathcal{H}$-class of $M(K\beta)$ divides some group $\mathcal{H}$-class of $M(K)$.
Clearly, a similar result can be deduced for $M(K\beta^{-1})$ and $M(K)$, by invoking Lemma~\ref{ab^-1:lema} instead of Lemma~\ref{ab:lema} throughout the proof of Lemma~\ref{lemma:Hclass:theta}. 

Since pseudovarieties are closed under division, the above discussion, together with Lemma~\ref{elem:monoid:same} and Lemma~\ref{compos:elem}, yields the following result.

\begin{thm}\label{thm:thebigone}
Let $\boldsymbol{\mathrm{H}}$ be a pseudovariety of groups. Then the pseudovariety of monoids $\overline{\boldsymbol{\mathrm{H}}}$ is stable under all automorphisms of $F_A$.
\end{thm}

Observe that if $\boldsymbol{\mathrm{H}}$ is the trivial pseudovariety (i.e. the pseudovariety of groups consisting of only trivial groups), then $\overline{\br{H}}$ is precisely the pseudovariety of aperiodic monoids, in view of \cite[Chapter 3 - Proposition 4.2]{pin}. And if $\boldsymbol{\mathrm{H}}$ is the pseudovariety of $p$-groups, then $\overline{\br{H}}$ is the pseudovariety of monoids all of whose subgroups are $p$-groups. Therefore, Theorem~\ref{thm:thebigone} generalizes the results obtained in \cite{inspo}, where it is implicitly shown that the pseudovarieties of aperiodic monoids and monoids all of whose subgroups are $p$-groups are stable under all automorphisms of $F_A$.

\subsection{Some characterizations}

We proceed by characterizing finitely generated subgroups $K \leq_{f.g.} F_A$ satisfying $M(K) \in \overline{\br{H}}$, for certain pseudovarieties of groups $\br{H}$.

Let $S = (k_n)_{n \geq 1}$ be a sequence of positive integers and denote by $\boldsymbol{\mathrm{V}}_S$ the pseudovariety of finite groups ultimately defined by the sequence $(x^{k_n}=1)_n$. In other words, a finite group $G$ belongs to $\boldsymbol{\mathrm{V}}_S$ if and only if there exists $p \geq 1$ such that $G$ satisfies the identity $x^{k_n}=1$ for all $n \geq p$. Observe that, since an identity of the type $x_n = y_n$ in a group is equivalent to $x_n y_n^{-1} = 1$, any pseudovariety defined by a sequence encompassing a single letter is of the type $\boldsymbol{\mathrm{V}}_S$ for some suitable $S$. Given two integers $m,n \geq 1$, we denote by $(m,n)$ the greatest common divisor of $m$ and $n$.

\begin{thm}\label{thm:psv:xkn}
Let $K \leq_{f.g.} F_A$ have Stallings automaton $\mathcal{S}(K) = (Q,A,\delta,q_0)$, and let $S = (k_n)_{n \geq 1}$ be a sequence of positive integers. The following conditions are equivalent:
\begin{enumerate}[label=(\roman*)]\setlength\itemsep{-0.1em}
\item $M(K) \in \overline{\boldsymbol{\mathrm{V}}}_S$.
\item $\exists p\geq 1, \forall q \in Q, \forall v \in \widetilde{A}^*, \forall n \geq 1, \forall m \geq p \quad q \cdot v^n = q \Rightarrow q \cdot v^{(n,k_m)} = q$.
\item $\exists p\geq 1, \forall x \in F_A, \forall n \geq 1, \forall m \geq p \quad x^n \in K \Rightarrow x^{(n,k_m)} \in K$.
\end{enumerate}
\end{thm}

\begin{proof}
(i) $\Rightarrow$ (ii): Assume that (i) holds. Let $q \in Q$, $v \in \widetilde{A}^*$ and $n \geq 1$ be such that $q \cdot v^n = q$. We define 
\begin{align*}
O(v) &= \bigcap_{m \geq 1} \dom \delta_{v^m} \\
\ell &= \min \{t \geq 1 \mid \forall o \in O(v) \; \, o \cdot v^t = o\},
\end{align*}
and choose $s \geq 1$ such that, for all $o \in \dom \delta_v \setminus O(v)$, $o \cdot v^s$ is not defined. Then we define $r = \ell s + 1$ and $w = v^r$. Since we can only read $w$ from vertices lying in $O(v)$, we deduce that $\delta_{w}$ belongs to a group $\mathcal{H}$-class of $M(K)$. Moreover, $q \in O(v)$ implies \[q \cdot w = q \cdot (v^\ell)^sv = q \cdot v,\] so $q \cdot w^i = q \cdot v^i$ for all $i \geq 1$. If $O(v) = \{o_1, \dots, o_m\}$ $(m \geq 1)$ and $d$ is the order of $\delta_w$, then $\delta_{w^d}$ is the restriction of the identity mapping fixing precisely $o_1, \dots, o_m$. Given $1 \leq i \leq m$, there is a smallest $d_i \geq 1$ such that $o_i \cdot v^{d_i} = o_i$, so we get disjoint cycles
\[\xymatrix{
o_i = o_{i,0} \ar[r]^-{v} & o_{i,1} \ar[r]^-{v} & \cdots \ar[r]^-{v} & o_{i,d_i-1} \ar@/^1pc/[lll]^-{v}
}\]
in $\mathcal{S}(K)$, where the vertices $o_{i,0}, o_{i,1}, \dots, o_{i,d_i-1} \in Q$ are all distinct. It is clear that $d$ is the least common multiple of $d_1, \dots, d_m$; in particular, $d_i \mid d$ for all $i$. Assuming, without loss of generality, that $q = o_1$, we must have $d_1 \mid n$, since $o_1 \cdot w^n = o_1 = o_1 \cdot w^{d_1}$. Furthermore, by hypothesis, there exists $p \geq 1$ such that ${\delta_w}^{k_m}$ is an idempotent for all $m \geq p$. This implies that $d_1 \mid k_m$ if $m \geq p$ because $d_1 \mid d$ and $d \mid k_m$. Thus, for all $m \geq p$, we obtain $d_1 \mid (n,k_m)$ and $q \cdot v^{(n,k_m)} = q \cdot w^{(n,k_m)} = q$, so (ii) follows.

(ii) $\Rightarrow$ (iii): Take $p$ satisfying \[\forall q \in Q, \forall v \in \widetilde{A}^*, \forall n \geq 1, \forall m \geq p \quad q \cdot v^n = q \Rightarrow q \cdot v^{(n,k_m)} = q.\] Let $x \in F_A$ and $n \geq 1$ be such that $x^n \in K$. We can assume that $x$ is reduced, seen as an element of $\widetilde{A}^*$, and write $x = uwu^{-1}$ for some $u,w \in \widetilde{A}$ with $w$ cyclically reduced. Then $q_0 \cdot \overline{x^n} = q_0 \cdot uw^nu^{-1} = q_0$. Letting $q = q_0 \cdot u$, it follows that $q \cdot w^n = q$. Hence, $q \cdot w^{(n,k_m)} = q$ for all $m \geq p$, by assumption. This entails \[ q_0 = q_0 \cdot uw^{(n,k_m)}u^{-1} = q_0 \cdot  \overline{x^{(n,k_m)}},\] so $x^{(n,k_m)} \in K$ for all $m \geq p$, and (iii) follows.

(iii) $\Rightarrow$ (ii): Take $p \geq 1$ satisfying \[\forall x \in F_A, \forall n \geq 1, \forall m \geq p \quad x^n \in K \Rightarrow x^{(n,k_m)} \in K.\] Let also $q \in Q$, $v \in \widetilde{A}^*$ and $n \geq 1$ be such that $q \cdot v^n = q$. Take a word $u \in \widetilde{A}^*$ satisfying $q_0 \cdot u = q$. Then $q_0 \cdot uv^nu^{-1} = q_0$ implies $(uvu^{-1})^n \in K$, and hence $(uvu^{-1})^{(n,k_m)} \in K$ for all $m \geq p$. 
The equalities $q \cdot u^{-1}\overline{uv^{(n,k_m)}u^{-1}}u = q$ and $\overline{u^{-1}\overline{uv^{(n,k_m)}u^{-1}}u} = \overline{v^{(n,k_m)}}$ allow us to conclude that $q \cdot v^{(n,k_m)} = q$ for all $m \geq p$, as we wanted.

(ii) $\Rightarrow$ (i): Take $p$ satisfying \[\forall q \in Q, \forall v \in \widetilde{A}^*, \forall n \geq 1, \forall m \geq p \quad q \cdot v^n = q \Rightarrow q \cdot v^{(n,k_m)} = q.\] Let $H = H_{\delta_v}$ be a group $\mathcal{H}$-class. If $d$ is the order of $\delta_v$ in $H$ then, for all $q \in \mathrm{dom}\, \delta_v$, we obtain $q \cdot v^d = q$. It follows that $q \cdot v^{(d,k_m)} = q$ for all $m \geq p$, by hypothesis. Since $(d,k_m) \leq d$, we obtain that $d = (d,k_m)$, and thus $d \mid k_m$ if $m \geq p$. We conclude that the order of every element in $H$ divides all the $k_m$ for $m \geq p$, thereby showing that $H \in \boldsymbol{\mathrm{V}}_S$. 
\end{proof}

Let $k \geq 1$ and denote by $\boldsymbol{\mathrm{B}}_k$ the \textbf{Burnside pseudovariety} given by $\boldsymbol{\mathrm{B}}_k = \llbracket x^k = 1 \rrbracket$. Considering the constant sequence $k_n = k$, we get the following corollary.

\begin{cor}\label{cor:psv:xk}
Let $K \leq_{f.g.} F_A$ have Stallings automaton $\mathcal{S}(K) = (Q,A,\delta,q_0)$. The following conditions are equivalent:
\begin{enumerate}[label=(\roman*)]\setlength\itemsep{-0.1em}
\item $M(K) \in \overline{\boldsymbol{\mathrm{B}}}_k$.
\item $\forall q \in Q, \forall v \in \widetilde{A}^*, \forall n \geq 1 \quad q \cdot v^n = q \Rightarrow q \cdot v^{(k,n)} = q$.
\item $\forall x \in F_A, \forall n \geq 1 \quad x^n \in K \Rightarrow x^{(k,n)} \in K$.
\end{enumerate}
\end{cor}

Let $\pi$ be a set of prime numbers and $\pi'$ its complement in the set of primes. A $\pi$-number is a number all of whose prime factors lie in $\pi$. Out of convenience, we will also consider $1$ to be a $\pi$-number, for every choice of $\pi$. A finite group $G$ is called a \textbf{$\pi$-group} if its order factors into primes from $\pi$. By Lagrange and Cauchy theorems, this is equivalent to saying that the order of every element in $G$ is a $\pi$-number. We denote by $\boldsymbol{\mathrm{G}}_\pi$ the pseudovariety of all finite $\pi$-groups. Observe that if $\pi = \{p_1, p_2,\dots\}$\footnote{If $\pi = \{p_1, p_2, \dots, p_m\}$ is finite, we consider $p_{m+1} = p_{m+2} = \dots = 1$.} and $S = \left((p_1\cdots p_n)^{n!}\right)_n$, then $\boldsymbol{\mathrm{G}}_\pi = \boldsymbol{\mathrm{V}}_S$ \cite[Proposition 7.1.16]{qtheory}. Hence, by Theorem~\ref{thm:psv:xkn}, we get a characterization of subgroups $K \leq_{f.g.} F_A$ satisfying $M(K) \in \overline{\boldsymbol{\mathrm{G}}}_\pi$. However, in this case, we present a somewhat simpler characterization, using essentially the same ideas. A proof can be found in \cite{tese}.

\begin{cor}
Let $K \leq_{f.g.} F_A$ with Stallings automaton $\mathcal{S}(K) = (Q,A,\delta,q_0)$. Let also $\pi$ be a set of prime numbers with complement $\pi'$ in the set of primes. Denote by $\Pi$ and $\Pi'$ the sets of $\pi$-numbers and $\pi'$-numbers, respectively. The following conditions are equivalent:
\begin{enumerate}[label=(\roman*)]\setlength\itemsep{-0.1em}
\item $M(K) \in \overline{\boldsymbol{\mathrm{G}}}_\pi$.
\item $\forall q \in Q, \forall v \in \widetilde{A}^*, \forall n \in \Pi' \quad q \cdot v^n = q \Rightarrow q \cdot v = q$.
\item $\forall x \in F_A, \forall n \in \Pi' \quad x^n \in K \Rightarrow x \in K$.
\end{enumerate}
\end{cor}

\section{Conjugacy conditions}\label{sec:conjug}

In this section, we discuss normal, malnormal and cyclonormal subgroups of $F_A$ using the transition monoid of the Stallings automaton.

\subsection{Normal subgroups}

We begin by analyzing the effect of conjugating a subgroup at the level of the Stallings automaton. Given a nontrivial subgroup $K \leq_{f.g.} F_A$, we denote by $\Gamma(K)$ the underlying labeled graph of $\St(K)$. In view of the abstract characterization of Stallings automata, $\Gamma(K)$ is of the form
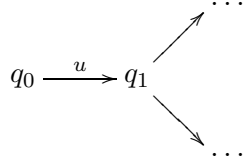
\vspace{2mm}
\[\xygraph{  
!{<0cm,0cm>;<1cm,0cm>:<0cm,1cm>::}
!{(0,0)}*+{q_0}="0"    
!{(1.5,0)}*+{q_1}="1" 
!{(2.5,1)}*+{\quad \dots}="." 
!{(2.5,-1)}*+{\quad \dots}=".."
"0":"1"^-u 
"1":"." "1":".."
} \]
where $q_1$ is either $q_0$ or the closest vertex to $q_0$ having outdegree (strictly) greater than $2$. We say that $q_0 \xrightarrow{\;u\;}$ is the \textbf{tail} of $\St(K)$ and the remaining graph is the \textbf{core} of $\St(K)$, which we denote by $\cl{C}(K)$. Note that the tail is not a graph, since the vertex $q_1$ does not belong to the tail, so it is empty if $q_1 = q_0$. Consequently, given $K \leq_{f.g.} F_A$, its Stallings automaton $\St(K)$ is either equal to $\cl{C}(K)$ with a vertex chosen to be the basepoint (if $K$ is normal in $F_A$ then this is the case); or it is obtained by ``gluing'' a tail $q_0 \xrightarrow{\;u\;}$ to $\cl{C}(K)$ and declaring $q_0$ to be the basepoint.

The next result states that two f.g. subgroups of a free group are conjugate if and only if their core graphs are isomorphic (for a proof, see \cite{kapovich}). 

\begin{thm}\label{thm:conjugacy}
Let $A$ be a finite alphabet and $H,K \leq_{f.g.} F_A$. Then $H$ and $K$ are conjugate subgroups if and only if $\cl{C}(H) \cong \cl{C}(K)$.
\end{thm}

\begin{exmp}
\normalfont Let $A = \{a,b,c\}$ and consider the Stallings automaton $\St(K)$
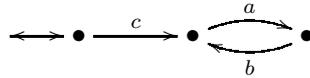
\[\xygraph{  
!{<0cm,0cm>;<1cm,0cm>:<0cm,1cm>::}
!{(-1,0)}*+{}="-1" 
!{(0,0)}*+{\bullet}="0"  
!{(1.5,0)}*+{\bullet}="1"    
!{(3,0)}*+{\bullet}="2" 
"-1":"0" "0":"-1"
"0":"1"^-c
"1":@/^/"2"^-{a} "2":@/^/"1"^-{b}
} \]
of a certain $K \leq_{f.g.} F_A$. Then its tail is
$\xygraph{  
!{<0cm,0cm>;<1cm,0cm>:<0cm,1cm>::}
!{(0,0)}*+{\bullet}="0"  
!{(1.5,0)}*+{}="1"    
"0":"1"^-c
}$
and we have

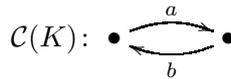
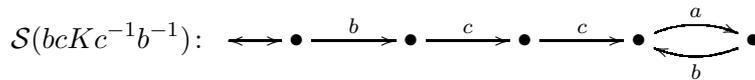
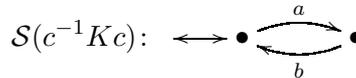
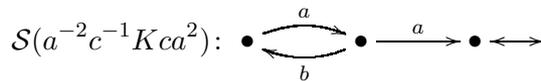
\[\cl{C}(K) \colon \xygraph{  
!{<0cm,0cm>;<1cm,0cm>:<0cm,1cm>::}
!{(0,0)}*+{\bullet}="1"    
!{(1.5,0)}*+{\bullet}="2" 
"1":@/^/"2"^-{a} "2":@/^/"1"^-{b}
} \]

\[\St(bcKc^{-1}b^{-1}) \colon \xygraph{  
!{<0cm,0cm>;<1cm,0cm>:<0cm,1cm>::}
!{(-4,0)}*+{}="-1" 
!{(-1.5,0)}*+{\bullet}="3" 
!{(-3,0)}*+{\bullet}="4" 
!{(0,0)}*+{\bullet}="0"  
!{(1.5,0)}*+{\bullet}="1"    
!{(3,0)}*+{\bullet}="2" 
"-1":"4" "4":"-1"
"0":"1"^-c
"1":@/^/"2"^-{a} "2":@/^/"1"^-{b}
"3":"0"^-c
"4":"3"^-b
} \]

\[\St(c^{-1}Kc) \colon \xygraph{  
!{<0cm,0cm>;<1cm,0cm>:<0cm,1cm>::}
!{(-1,0)}*+{}="-1" 
!{(0,0)}*+{\bullet}="1"    
!{(1.5,0)}*+{\bullet}="2" 
"-1":"1" "1":"-1"
"1":@/^/"2"^-{a} "2":@/^/"1"^-{b}
} \]

\[\St(a^{-2}c^{-1}Kca^2) \colon \xygraph{  
!{<0cm,0cm>;<1cm,0cm>:<0cm,1cm>::}
!{(4,0)}*+{}="-1" 
!{(0,0)}*+{\bullet}="1"    
!{(1.5,0)}*+{\bullet}="2" 
!{(3,0)}*+{\bullet}="3" 
"-1":"3" "3":"-1"
"1":@/^/"2"^-{a} "2":@/^/"1"^-{b}
"2":"3"^-a
} \]
\end{exmp}

\vspace{2mm}

We define an automorphism of a labeled directed graph to be a permutation $\sigma$ on its vertices such that, for every pair of vertices $p$ and $q$, if we have $p \xrightarrow{\; a \;} q$ then we also have $p\sigma \xrightarrow{\; a \;} q\sigma$. We say that a graph is \textbf{vertex-transitive} if, given two distinct vertices, there exists an automorphism of the graph taking one of the vertices to the other vertex. The following result is a consequence of Theorem~\ref{thm:conjugacy}. 

\begin{cor}\label{cor:vertextransitive}
Let $K$ be a nontrivial finitely generated subgroup of $F_A$. Then $K$ is normal in $F_A$ if and only if $K$ has finite index and $\cl{C}(K)$ is vertex-transitive.
\end{cor}

\begin{exmp}
\normalfont Let $A = \{a,b\}$. The automaton 
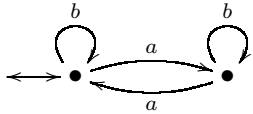
\[\xygraph{  
!{<0cm,0cm>;<1cm,0cm>:<0cm,1cm>::}
!{(-1,0)}*+{}="-1" 
!{(0,0)}*+{\bullet}="0"    
!{(2,0)}*+{\bullet}="1" 
"-1":"0" "0":"-1"
"0":@/^/"1"^{a} "1":@/^/"0"^{a}
"0":@(lu,ru)"0"^{b}
"1":@(lu,ru)"1"^{b}
} \]
is vertex-transitive and complete, so the corresponding subgroup is normal in $F_{A}$. However, if $K \leq_{f.g.} F_A$ is given by the Stallings automaton
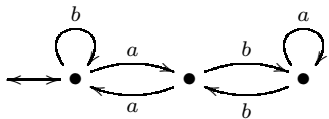
\[\xygraph{  
!{<0cm,0cm>;<1cm,0cm>:<0cm,1cm>::}
!{(-1,0)}*+{}="-1" 
!{(0,0)}*+{\bullet}="0"    
!{(1.5,0)}*+{\bullet}="1" 
!{(3,0)}*+{\bullet}="2" 
"-1":"0" "0":"-1"
"0":@/^/"1"^{a} "1":@/^/"0"^{a}
"0":@(lu,ru)"0"^{b}
"1":@/^/"2"^-{b}
"2":@/^/"1"^-{b}
"2":@(lu,ru)"2"^{a}
} \]
then $K$ is not a normal subgroup of $F_A$ because, despite being complete, $\St(K)$ is not vertex-transitive.
\end{exmp}

Before we characterize the normal subgroups of $F_A$ in terms of their transition monoids, we fix some notation and terminology regarding group actions. Let $X$ be a nonempty finite set and $S_X$ be the symmetric group over $X$. Given a group $G$, a \textbf{right action} of $G$ on $X$ is a function $\varphi \colon X \times G \rightarrow X$ satisfying $(x,1)\varphi = x$ and $(x,gh)\varphi = ((x,g)\varphi,h)\varphi$, for all $x \in X$ and $g,h \in G$. If $G \leq S_X$, then its elements are permutations on $X$, hence we get an action of $G$ on $X$ given by $(x,g) \mapsto xg$. For each $x \in X$, the \textbf{orbit} of $x$ is $xG = \{xg \mid g \in G \}$ and the \textbf{stabilizer} of $x$ is $G_x = \{g \in G \mid xg = x\} \leq G$. One easily checks that the orbits constitute a partition of $X$. We say that the action is \textbf{transitive} if it determines a single orbit. It is well known that $|xG| = [G:G_x]$, from which we derive the following result.

\begin{lem}\label{lemma:transitive}
Let $G$ be a finite group equipped with a transitive right action over a nonempty set $X$. The following conditions are equivalent:
\begin{enumerate}[label=(\roman*)]\setlength\itemsep{-0.1em}
\item $G_x = \{ 1\}$ for all $x \in X$.
\item $G_x = \{ 1\}$ for some $x \in X$.
\item $|G| = |X|$.
\end{enumerate}
\end{lem}

We are ready to state the aforementioned characterization of normal subgroups. 

\begin{thm}\label{thm:size:group}
Let $1 \neq K \leq_{f.g.} F_A$ and $M(K)$ be the transition monoid of its Stallings automaton $\mathcal{S}(K) = (Q,A,\delta,q_0)$. Then $K$ is normal in $F_A$ if and only if $M(K)$ is a group of size $|Q|$.
\end{thm}
\begin{proof}
Suppose first that $K \trianglelefteq F_A$. Since $K$ has finite index by Corollary~\ref{cor:vertextransitive}, $\mathcal{S}(K)$ is complete and $M(K)$ is a group, by Corollary~\ref{thm:completegroup}. Moreover, given $p,q \in Q$, there exists $u \in \widetilde{A}^*$ such that $p \cdot u = q$, as $\mathcal{S}(K)$ is connected, and hence $M(K) \leq S_Q$ acts transitively on $Q$. Let $\delta_u \in M(K)$ and $q \in Q$ be such that $q \cdot u = q$. Take an arbitrary vertex $p \in Q$ and consider an automorphism of $\mathcal{S}(K)$ inducing a bijection $\varphi \colon Q \rightarrow Q$ that maps $q$ to $p$. Then we have
\begin{align*}
p \cdot u = q\varphi \cdot u = (q \cdot u)\varphi = q\varphi = p,
\end{align*}
so $\delta_u$ fixes $p$ too. Hence, if $\delta_u \in M(K)$ is a transformation with a fixed point, it follows that $\delta_u = 1$. We conclude that $M(K)$ satisfies property (i) of Lemma~\ref{lemma:transitive}, thereby showing that $|M(K)| = |Q|$, by property (iii).

Assume now that $M(K) \leq S_Q$ is a group and $|M(K)| = |Q|$. As we remarked above, $M(K)$ acts transitively on $Q$ by connectedness of $\mathcal{S}(K)$ so, by Lemma~\ref{lemma:transitive}, the stabilizer of any state is trivial. Given any $p,q \in Q$, our goal is to define a bijection $\varphi \colon Q \rightarrow Q$ such that $p \varphi = q$ and, for all $r \in Q$ and $v \in \widetilde{A}$, the equality $r\varphi \cdot v = (r \cdot v)\varphi$ holds; this will show that $\cl{S}(K)$ is vertex-transitive. For every $r \in Q$, take $u_r \in \widetilde{A}^*$ satisfying $p \cdot u_r = r$ and define
\begin{align*}
\varphi \colon Q &\longrightarrow Q \\
r &\longmapsto q \cdot u_r.
\end{align*}
To check that $\varphi$ is well defined, observe that if $v \in \widetilde{A}^*$ is such that $p \cdot v = r$, then $p \cdot u_rv^{-1} = p$ and $\delta_{u_rv^{-1}} = 1$ by property (i) of Lemma~\ref{lemma:transitive}. Hence, $\delta_{u_r} = \delta_{v}$ and, in particular, $q \cdot u_r = q \cdot v$. Similarly one checks that $\varphi$ is an injection, for if $r = p \cdot u_r$ and $s = p \cdot u_s$ satisfy $r\varphi = s\varphi$, then $q \cdot u_r = q \cdot u_s$ implies $\delta_{u_r} = \delta_{u_s}$, and \[r = p \cdot u_r = p \cdot u_s = s.\] Since $Q$ is finite, we deduce that $\varphi$ is a bijection. Finally, note that for all $r \in Q$ and $v \in \widetilde{A}$, we have
\begin{align*}
r\varphi \cdot v &= (q \cdot u_r) \cdot v \\
&= q \cdot u_rv \\
&= (p \cdot u_rv)\varphi \\
&= ((p \cdot u_r)\cdot v)\varphi \\
&= (r \cdot v)\varphi.
\end{align*}
The desired conclusion now follows.
\end{proof}

\subsection{Malnormal subgroups}

Given any group $G$ and $K \leq G$, we say that $K$ is \textbf{malnormal} in $G$ if $gKg^{-1} \cap K = 1$ for all $g \in G \setminus K$. In the case of free groups, there is the following well-known result \cite[Theorem~9.10]{kapovich}.

\begin{lem}\label{thm:malnormal:stallings}
A subgroup $K \leq_{f.g.} F_A$ is malnormal if and only if, in its Stallings automaton $\mathcal{S}(K) = (Q,A,\delta,q_0)$, there do not exist two distinct vertices $q, q' \in Q$ and a nonempty reduced word $u$ such that $u$ labels a loop at both $q$ and $q'$.
\end{lem}

Before stating our result regarding malnormality, we present some definitions and make a few observations. Given an inverse monoid $M$, the \textbf{natural partial order} on $M$ is the partial order $\leq$ given by
\[x \leq y \Leftrightarrow x = ey \; \text{for some} \; e \in E(M).\] If the elements of $M$ are injective partial transformations, then $f,g \in M$ satisfy $f \leq g$ if and only if $g|_{\mathrm{dom}\, f} = f$. Note also that a partial transformation is an idempotent if and only if it is a restriction of the identity. The next result is just another simple remark.

\begin{lem}\label{lema:idemp}
Let $K \leq_{f.g.} F_A$ and $x \in \widetilde{A}^*$ be such that $\delta_x \in E(M(K))$. Then, for all $r \in \widetilde{A}^*$ and $n \geq 1$, the condition $\delta_{rx^nr^{-1}} = \delta_{rxr^{-1}}\in E(M(K))$ holds.
\end{lem}

Given a partially ordered set $(X, \leq)$ and $Y = \{y_1, \dots, y_n \} \subseteq X$ for some $n \geq 1$, we say that $Y$ is a \textbf{chain} on $X$ if there exists $\sigma \in S_n$ satisfying $y_{1\sigma} < \cdots < y_{n\sigma}$. We say that a chain of the type $y_1 < y_2 < \dots < y_k$ ($k \geq 1$) has length $k$. 

Finally, a word $u \in \widetilde{A}^*$ is called \textbf{cyclically reduced} if $uu$ is a reduced word, which is equivalent to saying that $u$ is reduced and there exist no $a \in \widetilde{A}$ and $v \in \widetilde{A}^*$ satisfying $u = ava^{-1}$. It is easy to see that every reduced word $u$ can be uniquely written as $u = xvx^{-1}$ for some reduced word $x$ and cyclically reduced word $v$. Observe that if a word $uw$ is cyclically reduced, then $wu$ is cyclically reduced too.

Now we characterize the malnormal subgroups of $F_A$ in terms of their transition monoids.

\begin{thm}\label{thm:malnormal:me}
Given a nontrivial subgroup $ K <_{f.g.} F_A$, consider the monoid homomor\-phism 
\begin{align*}
\Delta \colon {\widetilde{A}}^* &\longrightarrow M(K) \\
u &\longmapsto \delta_u
\end{align*}
and let $R_A$ be the set of all reduced words over $\widetilde{A}$. Denote by $E$ the set of idempotents of $(R_A \setminus \{1\}) \Delta$. Consider the restriction of the natural partial order on $M(K)$ to $E$, and let $k$ be the size of a maximal chain on $E$. Then $K$ is malnormal if and only if $k=2$ and $|E| = |Q| + 1$.
\end{thm}

\begin{proof} 
First observe that, by construction of $\mathcal{S}(K) = (Q,A,\delta,q_0)$, there exist a nonempty cyclically reduced word $u \in R_A$ and $q \in Q$ such that $u$ labels a loop at $q$. Hence, such $\delta_u$ is not the empty transformation on $Q$. Since $M(K)$ is finite, there exists $n \geq 1$ such that $\delta_u^n \in E$. So if $k=1$, then $0 \not \in E$ and, consequently, $0 \not \in M(K)$. This implies that $\mathcal{S}(K)$ is a complete automaton by Lemma \ref{lema:complete}, so $M(K)$ is a finite group. As $K \neq F_A$, it follows that $\mathcal{S}(K)$ has at least two vertices. Given $a \in A$, there exists some $n \geq 1$ such that $\delta_a^n = \rr{id}_Q$, and clearly $a^n \in R_A \setminus \{1\}$. Hence, if $q$ and $q'$ are two distinct vertices in $\mathcal{S}(K)$, the word $a^n$ labels a loop at $q$ and labels a loop at $q'$. This entails that $K$ is not malnormal.

Now suppose that $k \geq 3$. Then there exist at least three distinct elements $\delta_u, \delta_v, \delta_w \in E$ satisfying $\delta_u < \delta_v < \delta_w$. Accordingly, $\delta_w$ fixes at least two distinct vertices, say $q$ and $q'$, which implies that $w \in R_A \setminus \{1\}$ labels a loop at both $q$ and $q'$. Hence, $K$ is not malnormal.

All in all, if $K$ is malnormal, then $k=2$. Note that, in this case, $0 \in E$; otherwise, by Lemma~\ref{lema:complete}, $M(K)$ would be a group, and thus it would have a unique idempotent. So there exist distinct nonempty sets $Q_1, Q_2, \dots, Q_\ell \subseteq Q$ ($\ell \geq 1$) such that \[E = \{0, \mathrm{id}_{Q_1}, \dots, \mathrm{id}_{Q_\ell}\},\] yielding $|E| = \ell + 1$. Now, given any $q \in Q$, by construction of $\mathcal{S}(K)$ there exists some $u \in R_A \setminus \{1\}$ satisfying $q \cdot u = q$, which means that $\delta_u$ fixes $q$. As $\mathcal{S}(K)$ is finite, there exists $n \geq 1$ for which $\overline{u^n} \neq 1$ and $\delta_{\overline{u^n}}$ is an idempotent fixing $q$. This shows that $q \in Q_i$ for some $1 \leq i \leq \ell$, and therefore \[Q = Q_1 \cup Q_2 \cup \dots \cup Q_\ell.\]
If $|Q|>1$, malnormality of $K$ implies that $|Q_i| = 1$, for $1 \leq i \leq \ell$, which gives $|Q| = \ell$ and $|E| = |Q| + 1$; if $|Q| = 1$, then $E = \{0, \mathrm{id}_Q \}$ and $|E| = |Q|+1$ too.

Finally, keeping the notation introduced in the preceding paragraph, suppose that $k=2$ and $|E| = |Q| + 1$ (so $\ell = |Q|$). In order to prove that $K$ is malnormal, we only need to show that, for all $1 \leq i,j \leq \ell$, \[ i \neq j  \Rightarrow Q_i \cap Q_j = \emptyset. \]
Indeed, if that is the case, we get \[|Q| = |Q_1| + |Q_2| + \cdots + |Q_\ell| \geq \ell.\] So, under the hypothesis that $|Q| = \ell$, all the $Q_i$ must be singletons, as they are nonempty. Hence, if $K$ is not malnormal, then there exist distinct states $q, q' \in Q$ and $u \in R_A \setminus \{1\}$ such that $q \cdot u = q$ and $q' \cdot u = q'$. Also, there exists $n \geq 1$ such that $\delta_{\overline{u^n}} = \mathrm{id}_{Q_j} \in E$ for some $1 \leq j \leq \ell$. Since $q,q' \in Q_j$, it follows that $|Q_j| > 1$, which contradicts our assumption. Therefore, if all the $Q_i$ ($1 \leq i \leq \ell)$ are disjoint, then $K$ is a malnormal subgroup of $F_A$.

We proceed by proving the disjointness claim by contradiction. Suppose that there exist $1 \leq i < j \leq \ell$ and $u,v \in R_A \setminus \{1\}$ such that $\delta_u = \mathrm{id}_{Q_i}$, $\delta_v = \mathrm{id}_{Q_j}$ and $q \in Q_i \cap Q_j$. Note that, by assumption, $Q_i \neq Q_j$. As before, we can write $u = \widetilde{u} x \widetilde{u}^{-1}$ and $v = \widetilde{v} y \widetilde{v}^{-1}$, for some $\widetilde{u}, \widetilde{v} \in R_A$ and $x, y \in R_A \setminus \{1\}$, so that both $x$ and $y$ are cyclically reduced words. Moreover, let $t \in R_A$ be the longest prefix shared by $\widetilde{u}$ and $\widetilde{v}$ and write $\widetilde{u} = tr$, $\widetilde{v} = ts$ for some appropriate $r, s \in R_A$. Then $u = trxr^{-1}t^{-1}$ and $v = tsys^{-1}t^{-1}$. 

Since $M(K)$ is finite, there exists $n \geq 1$ such that $\delta_{x^n}, \delta_{y^n} \in E$, which implies \[\delta_{rx^nr^{-1}}, \delta_{sy^ns^{-1}} \in E,\] by Lemma~\ref{lema:idemp}. Letting $w = trx^nr^{-1}t^{-1}$ and $z = tsy^ns^{-1}t^{-1}$, we also have $\delta_w = \delta_u$ and $\delta_z = \delta_v$. Hence, we deduce that $\delta_w \neq \delta_z$ and $\delta_{rx^nr^{-1}} \neq  \delta_{sy^ns^{-1}}$.

We now analyse three possible cases:
\begin{itemize}
\item \textbf{Case 1:} $r=1$ and $s=1$.

We have $u = txt^{-1}$ and $v = tyt^{-1}$, with $x \neq y$ nonempty cyclically reduced words. If $xy \not \in R_A$, then we can write $x = \widetilde{x}a$ and $y = a^{-1}\widetilde{y}$ for some $a \in \widetilde{A}$ and $\widetilde{x}, \widetilde{y} \in R_A$. In that case, the last letter of $y$ cannot be $a$, and thus $xy^{-1} \in R_A \setminus \{1\}$. So there exists $\varepsilon \in \{-1,1\}$ such that $xy^\varepsilon \in R_A \setminus \{1\}$. However, this yields a contradiction, for we would get a chain
\begin{gather*}
0 < \delta_{x^ny^{\varepsilon n}} < \delta_{x^n} \\
\text{or} \\
0 < \delta_{x^ny^{\varepsilon n}} < \delta_{y^{\varepsilon n}}
\end{gather*}
in $E$, which violates the assumption that $k=2$. Indeed, $\delta_{x^ny^{\varepsilon n}} \neq 0$ because $x^ny^{\varepsilon n}$ labels a loop at $(q,t)\delta$, as $x$ and $y$ do so; and $\delta_{y^{\varepsilon n}} \neq \delta_{x^n}$ entails $\delta_{x^n y^{\varepsilon n}} < \delta_{x^n}$ or $\delta_{x^n y^{\varepsilon n}} < \delta_{y^{\varepsilon n}}$.

\item \textbf{Case 2:} $r = 1$ and $s \neq 1$.

We have $u = txt^{-1}$ and $v = tsys^{-1}t^{-1}$, with $x \neq y$ nonempty cyclically reduced words. If $xsys^{-1} \not \in R_A$, then we can write $x = \widetilde{x}a$ and $s = a^{-1}\widetilde{s}$ for some $a \in \widetilde{A}$ and $\widetilde{x}, \widetilde{s} \in R_A$. In that case, the first letter of $x$ cannot be $a^{-1}$, and thus $sys^{-1}x \in R_A \setminus \{1\}$. Either way, we can use an argument similar to the above to get a chain
\begin{gather*}
0 < \delta_{x^n sy^ns^{-1}} < \delta_{x^n} \;\; \text{or} \;\; 0 < \delta_{x^n sy^ns^{-1}} < \delta_{sy^ns^{-1}} \\
\text{or} \quad \; \; \\
0 < \delta_{sy^ns^{-1}x^n } < \delta_{x^n} \;\; \text{or} \;\; 0 < \delta_{sy^ns^{-1}x^n } < \delta_{sy^ns^{-1}}
\end{gather*}
in $E$, yielding a contradiction again. The case $r \neq 1$ and $s=1$ is analogous.

\item \textbf{Case 3:} $r\neq1$ and $s \neq 1$.

We have $u = trxr^{-1}t^{-1}$ and $v = tsys^{-1}t^{-1}$, with $x \neq y$ nonempty cyclically reduced words. If $rxr^{-1}sys^{-1} \not \in R_A$, then we can write $r = a\widetilde{r}$ and $s = a\widetilde{s}$ for some $a \in \widetilde{A}$ and $\widetilde{r}, \widetilde{s} \in R_A$. However, that is not possible, for otherwise $ta$ would be a longer prefix common to $\widetilde{u}$ and $\widetilde{v}$. So $rxr^{-1}sys^{-1} \in R_A$ and we get a chain
\begin{gather*}
0 < \delta_{rx^nr^{-1} sy^ns^{-1}} < \delta_{rx^nr^{-1}} \\
\text{or} \\
0 < \delta_{rx^nr^{-1} sy^ns^{-1}} < \delta_{sy^ns^{-1}}
\end{gather*}
in $E$, contradicting the assumption that $k=2$.
\end{itemize}

In any case, we always arrive at a contradiction, which means that if $1 \leq i < j \leq \ell$, then $Q_i \cap Q_j = \emptyset$. This, as we remarked, concludes the proof of the theorem.
\end{proof}

\subsection{Cyclonormal subgroups}

Given any group $G$ and $K \leq G$, we say that $K$ is \textbf{cyclonormal} if $gKg^{-1} \cap K$ is a cyclic group for every $g \in G \setminus K$. 

Let $K \leq_{f.g.} F_A$ with $\St(K) = (Q,A,\delta,q_0)$, and let $p, q \in Q$. Denoting by $\phi \colon \widetilde{A}^* \rightarrow F_A$ the quotient morphism, we define \[L_{(p,q)} = \{u\phi \mid u \in \widetilde{A}^* \text{ labels a loop at } (p,q) \text{ in } \St(K) \times \St(K) \} \leq F_A.\] 
Here, $\St(K) \times \St(K)$ denotes the \textbf{direct product} of $\St(K)$ with itself, so its vertex set is $Q \times Q$ and there exists an edge $(p,p') \xrightarrow{\; a \;} (q,q')$ if and only if we have edges $p \xrightarrow{\; a \;} q$ and $p' \xrightarrow{\; a \;} q'$ in $\St(K)$. In Stallings' terminology, this corresponds to the pull-back of $\St(K)$ with itself.

It is clear that if $|Q| = 1$, then $K$ is cyclonormal; and if $|A| = 1$, then every subgroup of $F_A$ is cyclonormal. Below we present a result regarding cyclonormality in free groups.

\begin{thm}\rm{\cite{kapovich}}\label{thm:cyclo:stallings}
Let $K \leq_{f.g.} F_A$ be a subgroup of $F_A$ and $\mathcal{S}(K) = (Q,A,\delta,q_0)$. The following conditions are equivalent:
\begin{enumerate}[label=(\roman*)]\setlength\itemsep{-0.1em}
\item $K$ is cyclonormal.
\item For each pair of distinct vertices $(p,q) \in Q \times Q$, the subgroup $L_{(p,q)}$ is cyclic.
\item For every $p \neq q \in Q$, there exists $u \in R_A$ such that $p \cdot u = p$, $q \cdot u = q$, and for all $v \in R_A$ satisfying $p \cdot v = p$ and $q \cdot v = q$, the equality $v = \overline{u^m}$ holds for some $m \in \mathbb{Z}$.
\end{enumerate}
\end{thm}

\color{black}

The following lemma is straightforward but contains a useful observation.

\begin{lem}\label{lemma:combinatorics}
Let $A$ be a set of size $n>2$ and $A_1, \dots, A_\ell \subseteq A$ ($\ell \geq 1$) be distinct nonempty subsets such that $\cup_{i=1}^\ell A_i = A$. If there are more than $\binom{n}{2}$ nonsingletons among the $A_i$, then there exist distinct $1 \leq i,j \leq \ell$ satisfying $|A_i \cap A_j| > 1$.
\end{lem}

We are ready to prove our main result.

\begin{thm}\label{thm:cyclonormal:me}
Suppose that $|A| \geq 2$ and let $1 \neq K <_{f.g.} F_A$ be a cyclonormal subgroup of $F_A$ with Stallings automaton $\mathcal{S}(K) = (Q,A,\delta,q_0)$. Let $E$ and $k\geq 1$ be as in Theorem~\ref{thm:malnormal:me}. Then either $k=2$ or $k=3$. Moreover, if $|Q| = n > 2$, then $k=2$ implies $|E| \leq {n \choose 2} + 1 $; and if $k=3$, then $|E| \leq n + {n \choose 2} + 1$. 
\end{thm}
\begin{proof}
If $k=1$, then $M(K)$ is a finite group and $\mathcal{S}(K)$ is a complete automaton, as we observed in the proof of Theorem~\ref{thm:malnormal:me}. Consequently, given two distinct letters $a,b \in A$, there exist $m,n \in \mathbb{Z}$ such that, given any pair of vertices $p$ and $q$, the words $a^m$ and $b^n$ label loops at both $p$ and $q$. However, there are no $u \in R_A$ and $s,t \in \mathbb{Z}$ satisfying $a^m = \overline{u^s}$ and $b^n = \overline{u^t}$, as the first letter of $u$ would have to be simultaneously equal to $a$ and $b$. Hence, if $K$ is cyclonormal, then $k > 1$.

Now let $Q_1, Q_2, \dots, Q_\ell \subseteq Q$ ($\ell \geq 1$) be distinct nonempty sets satisfying \[E = \{0, \mathrm{id}_{Q_1}, \dots, \mathrm{id}_{Q_\ell}\}.\] Since $K$ is cyclonormal, there cannot exist $i \neq j$ such that $|Q_i \cap Q_j| \geq 2$. Otherwise, we could find $\delta_v \neq \delta_w \in E$ fixing at least two vertices $p, q \in Q$, and we may assume that $\mathrm{dom}\, \delta_w \not\subset \mathrm{dom}\, \delta_v$. Suppose that there were $u \in R_A$ labelling a loop at both $p$ and $q$ and $s,t \in \mathbb{Z}$ satisfying $v = \overline{u^s}$ and $w = \overline{u^t}$. Then taking $r \in \mathrm{dom}\, \delta_w \setminus \mathrm{dom}\, \delta_v$, we would get $r \cdot u^{|t|} = r \cdot u^t = r \cdot w = r$, whence $r \cdot u^{|t|s} = r$. However, we would also get $r \cdot \overline{u^s} = r \cdot v \not \in Q$, which implies $r \cdot u^s \not \in Q$ and $r \cdot u^{|t|s} \not \in Q$, a contradiction. This allows us to conclude that $k \leq 3$, for otherwise there would exist $\delta_v < \delta_w$ whose domains would have at least two vertices in common.

Moreover, if $|Q| = n > 2$, $k=2$ and $|E| > {n \choose 2} + 1 $, then $\ell > {n \choose 2}$. If there are no singletons among the $Q_i$, then Lemma~\ref{lemma:combinatorics} holds for $A = Q$ and $A_i = Q_i$ ($1 \leq i \leq \ell$); therefore, there exist distinct $i$ and $j$ such that $|Q_i \cap Q_j| \geq 2$, a contradiction since $K$ is cyclonormal. Otherwise, assume that $Q_1, \dots, Q_t$ ($1 \leq t \leq \ell$) are singletons, say $Q_1 = \{q_1\}, \cdots, Q_t = \{q_t\}$, for some $q_1, \dots, q_t \in Q$. Let $\sigma \in S_Q$ be a permutation on $Q$ whose cycles have length greater than $2$ (this is possible since $|Q| > 2$). For $1 \leq i \leq t$, let $B_i = \{q_i, q_i\sigma\}$. Then the nonsingular sets $B_1,\dots,B_t,Q_{t+1},\dots,Q_\ell$ cover $Q$, the $B_i$ are all distinct by the choice of $\sigma$, and no $B_i$ is equal to a $Q_j$ if $j \neq i$, as $k=2$ implies $Q_i \not\subset Q_j$. By Lemma~\ref{lemma:combinatorics}, we conclude that there exist $1 \leq i \leq t$ and $t+1 \leq j \leq \ell$ for which $|B_i \cap Q_j| \geq 2$, thereby implying $Q_i \subset Q_j$, which is absurd. The contradiction we arrived at followed from the assumption that $|E| > \binom{n}{2}$. We conclude that $|E| \leq {n \choose 2} + 1 $ as claimed.

Finally, if $|Q|=n>2$, $k=3$ and $\ell > n + \binom{n}{2}$, there are certainly more than $\binom{n}{2}$ nonsingletons among the $Q_i$. It follows from the previous lemma that there exist distinct $1 \leq i,j \leq \ell$ such that $|Q_i \cap Q_j| \geq 2$, which contradicts the fact that $K$ is cyclonormal. We conclude that, in this case, $\ell \leq n + \binom{n}{2}$, and thus $|E| \leq n+{n \choose 2}+1$.
\end{proof} 

The following example shows that $M(K)$ and $E$ alone are not sufficient to completely characterize cyclonormality. Let $A = \{a,b,c\}$ and consider the subgroups $H,K \leq_{f.g.} F_A$ with Stallings automata
\[\mathcal{S}(H): \xygraph{  
!{<0cm,0cm>;<1cm,0cm>:<0cm,1cm>::}  
!{(-1,0)}*+{}="-1" 
!{(0,0)}*+{1}="1"    
!{(2,0)}*+{2}="2" 
"-1":"1" "1":"-1"  
"1":"2"^-{c}
"1":@(lu,ru)"1"^{a}
"2":@(lu,ru)"2"^{a}
} \] 
\[\mathcal{S}(K): \xygraph{  
!{<0cm,0cm>;<1cm,0cm>:<0cm,1cm>::}  
!{(-1,0)}*+{}="-1" 
!{(0,0)}*+{1}="1"    
!{(2,0)}*+{2}="2" 
"-1":"1" "1":"-1"  
"1":"2"^-{c}
"1":@(lu,ru)"1"^{a,b}
"2":@(lu,ru)"2"^{a,b}
} \]  
It is clear that $M(K) = M(H)$ and the sets $E$ (as in Theorem~\ref{thm:malnormal:me}) are also equal. However, $H$ is cyclonormal whereas $K$ is not.

\color{black}

\section{Future research}

Regarding future directions of work, it would be interesting to determine more properties $P$ (not necessarily related to pseudovarieties) satisfying the following condition: given a subgroup $K \leq_{f.g.} F_A$, if $M(K)$ satisfies property $P$, then $M(K\varphi)$ also satisfies property $P$, for all $\varphi \in \mathrm{Aut}(F_A)$. The next step would be to link the algebraic property $P$ of $M(K)$ to an algebraic property of $K$ as a subgroup of $F_A$.

Moreover, it could also be interesting to study the effect of various operators acting on the lattice of f.g. subgroups of a free group at the level of the transition monoids of Stallings automata.

Finally, one could pursue these same goals regarding structures that generalize Stallings automata for some wider classes of groups \cite{next5,next1,next4,next2,next3}.

\subsection*{Acknowledgments}
The author would like to thank the anonymous referees for providing helpful comments leading to the improvement of the paper. A special thanks to Professor Pedro Silva for all the encouragement and fruitful discussions.
The author was partially supported by the grant UIBD/MAT/00144/2020 (Ref. I) from CMUP, which is financed by national funds through FCT -- Funda\c c\~ao para a Ci\^encia e a Tecnologia, I.P.

\end{document}